\documentclass[11pt, reqno]{amsart}

\usepackage{amssymb,amsthm}
\usepackage[T1]{fontenc}
\usepackage{pdfpages}

\usepackage[colorlinks=True]{hyperref}

\evensidemargin\oddsidemargin

\begin{document}

\baselineskip=16pt
\setcounter{page}{1}
    
\newtheorem{Theorem}{Theorem}[section]
\newtheorem{Lemma}[Theorem]{Lemma}
\newtheorem{Cor}[Theorem]{Corollary}
\newtheorem{Prop}[Theorem]{Proposition}
\newtheorem{Rem}[Theorem]{Remark}
\newtheorem{Exa}[Theorem]{Example}
\newtheorem{Def}[Theorem]{Definition}


\newcommand{\ds}[1]{{\displaystyle{#1 }}}

\renewcommand{\leq}{\leqslant}
\renewcommand{\geq}{\geqslant}

\newcommand{\ind}{\mathds{1}}

\newcommand{\dd}{{\mathrm d}}
\newcommand{\ee}{{\mathrm e}}
\newcommand{\Span}{\mathrm{Span}}
\newcommand{\Supp}{\mathrm{Supp}}

\newcommand{\rc}{{\mathrm c}}

\newcommand{\tb}{{\mathrm b}}

\newcommand{\dC}{\mathbb{C}}
\newcommand{\dE}{\mathbb{E}}
\newcommand{\dN}{\mathbb{N}}
\newcommand{\dP}{\mathbb{P}}
\newcommand{\dR}{\mathbb{R}}
\newcommand{\dV}{\mathbb{V}}
\newcommand{\dZ}{\mathbb{Z}}

\newcommand{\cA}{\mathcal{A}}
\newcommand{\cB}{\mathcal{B}}
\newcommand{\cC}{\mathcal{C}}
\newcommand{\cF}{\mathcal{F}}
\newcommand{\cG}{\mathcal{G}}
\newcommand{\cI}{\mathcal{I}}
\newcommand{\cL}{\mathcal{L}}
\newcommand{\cM}{\mathcal{M}}
\newcommand{\cN}{\mathcal{N}}
\newcommand{\cP}{\mathcal{P}}
\newcommand{\cQ}{\mathcal{Q}}
\newcommand{\cR}{\mathcal{R}}
\newcommand{\cU}{\mathcal{U}}
\newcommand{\cW}{\mathcal{W}}

\newcommand{\bD}{\mathbf{D}}
\newcommand{\bL}{\mathbf{L}}
\newcommand{\bM}{\mathbf{M}}
\newcommand{\bS}{\mathbf{S}}
\newcommand{\bU}{\mathbf{U}}
\newcommand{\bW}{\mathbf{W}}
\newcommand{\bX}{\mathbf{X}}
\newcommand{\bY}{\mathbf{Y}}
\newcommand{\bZ}{\mathbf{Z}}
\newcommand{\be}{\mathbf{e}}
\newcommand{\bbf}{\mathbf{f}}
\newcommand{\bm}{\mathbf{m}}
\newcommand{\bu}{\mathbf{u}}
\newcommand{\bv}{\mathbf{v}}
\newcommand{\bw}{\mathbf{w}}
\newcommand{\bx}{\mathbf{x}}
\newcommand{\by}{\mathbf{y}}
\newcommand{\bz}{\mathbf{z}}

\newcommand{\PENT}[1]{{\lfloor#1\rfloor}} 
\newcommand{\PFRAC}[1]{{\lceil#1\rceil}} 
\newcommand{\ABS}[1]{{{\left| #1 \right|}}} 
\newcommand{\BRA}[1]{{{\left\{#1\right\}}}} 
\newcommand{\SCA}[1]{{{\left<#1\right>}}} 
\newcommand{\NRM}[1]{{{\left\| #1\right\|}}} 
\newcommand{\PAR}[1]{{{\left(#1\right)}}} 
\newcommand{\BPAR}[1]{{{\biggl(#1\biggr)}}} 
\newcommand{\BABS}[1]{{{\biggl|#1\biggr|}}} 
\newcommand{\SBRA}[1]{{{\left[#1\right]}}} 
\newcommand{\DSBRA}[1]{{{\llbracket#1\rrbracket}}} 
\newcommand{\VT}[1]{{{\| #1\|}_{\mbox{{\scriptsize VT}}}}} 
\newcommand{\LIP}[1]{{\|#1\|_{\mathrm{Lip}}}} 

\def\a{\alpha}
\def\b{\beta}
\def\B{{\bf B}} 
\def\C{{\bf C}} 
\def\CC{{\mathbb{C}}}
\def\cG{{\mathcal{G}}} 
\def\cH{{\mathcal{H}}} 
\def\cI{{\mathcal{I}}} 
\def\cS{{\mathcal{S}}}
\def\UU{{\mathcal{U}}}
\def\ca{c_{\a}}
\def\ka{\kappa_{\a}}
\def\coa{c_{\a, 0}}
\def\cua{c_{\a, u}}
\def\cL{{\mathcal{L}}} 
\def\cM{{\mathcal{M}}} 
\def\Ea{E_\a}
\def\eps{{\varepsilon}} 
\def\esp{{\mathbb{E}}} 
\def\Ga{{\Gamma}} 
\def\GG{{\bf \Gamma}} 
\def\e{{\rm e}}
\def\ii{{\rm i}}
\def\L{{\bf L}}
\def\lbd{\lambda}
\def\lacc{\left\{}
\def\lcr{\left[}
\def\lpa{\left(}
\def\lva{\left|}
\def\M{{\bf M}}
\def\NN{{\mathbb{N}}} 
\def\pb{{\mathbb{P}}}
\def\QQ{{\mathbb{Q}}} 
\def\R{{\bf R}}
\def\rl{{\mathbb{R}}}
\def\racc{\right\}}
\def\rpa{\right)}
\def\rcr{\right]}
\def\rva{\right|}
\def\W{{\bf W}}
\def\X{{\bf X}}
\def\XX{{\mathcal X}}
\def\YY{{\mathcal Y}}
\def\Y{{\bf Y}}
\def\V{{\bf V}_\a}
\def\Un{{\bf 1}}
\def\S{{\bf S}}
\def\A{{\bf A}}
\def\G{{\bf G}}
\def\AA{{\mathcal A}}
\def\hAA{{\widehat \AA}}
\def\hL{{\widehat L}}
\def\T{{\bf T}}

\def\claw{\stackrel{d}{\longrightarrow}}
\def\elaw{\stackrel{d}{=}}
\def\qed{\hfill$\square$}

\newcommand*\pFqskip{8mu}
\catcode`,\active
\newcommand*\pFq{\begingroup
        \catcode`\,\active
        \def ,{\mskip\pFqskip\relax}%
        \dopFq
}
\catcode`\,12
\def\dopFq#1#2#3#4#5{%
        {}_{#1}F_{#2}\biggl[\genfrac..{0pt}{}{#3}{#4};#5\biggr]%
        \endgroup
}

\newcommand{\helene}[1]{\textcolor{purple}{#1}}
\newcommand{\kilian}[1]{\textcolor{brown}{#1}}


\title[Superdiffusive limit of the elephant random walk]{On the limit law of the superdiffusive elephant~random~walk}

\date{\today}
\author[H.~Guérin]{Hélène Guérin}

\address{Université du Québec à Montréal, Département de mathématiques, Montréal, Canada. {\em Email}: {\tt guerin.helene@uqam.ca}}

\author[L.~Laulin]{Lucile Laulin}

\address{Université Paris-Nanterre, Modal'X, Nanterre, France. {\em Email}: {\tt lucile.laulin@math.cnrs.fr}}

\author[K.~Raschel]{Kilian Raschel} 

\address{CNRS, International Research Laboratory France-Vietnam in mathematics and its applications, Vietnam Institute for Advanced Study in Mathematics, Hanoï, Vietnam. {\em Email}: {\tt raschel@math.cnrs.fr}}

\author[T.~Simon]{Thomas Simon}

\address{Universit\'e de Lille, Laboratoire Paul Painlev\'e, Lille, France. {\em Email}: {\tt thomas.simon@univ-lille.fr}}



\begin{abstract} 
When the memory parameter of the elephant random walk is above a critical threshold, the process becomes superdiffusive and, once suitably normalised, converges to a non-Gaussian random variable. In a recent paper by the three first authors, it was shown that this limit variable has a density and that the associated moments satisfy a nonlinear recurrence relation. In this work, we exploit this recurrence to derive an asymptotic expansion of the moments and the asymptotic behaviour of the density at infinity. In particular, we show that an asymmetry in the distribution of the first step of the random walk leads to an asymmetry of the tails of the limit variable. These results follow from a new, explicit expression of the Stieltjes transformation of the moments in terms of special functions such as hypergeometric series and incomplete beta integrals. We also obtain other results about the random variable, such as unimodality and, for certain values of the memory parameter, log-concavity.
\end{abstract}

\maketitle

\section{Introduction and statement of the results}
\label{sec:introduction}

The one-dimensional elephant random walk (ERW)\ was introduced in 2004 by Schütz and Trimper \cite{Schutz2004} to see how memory could induce subdiffusion in random walk processes. It turned out that the ERW is always at least diffusive; nevertheless, the simple definition of the process, together with the underlying depth of the results, has led to a great deal of interest from mathematicians over the last two decades. 

\begin{figure}[ht!]
    \centering
    \includegraphics[width=7cm,height=5.5cm]{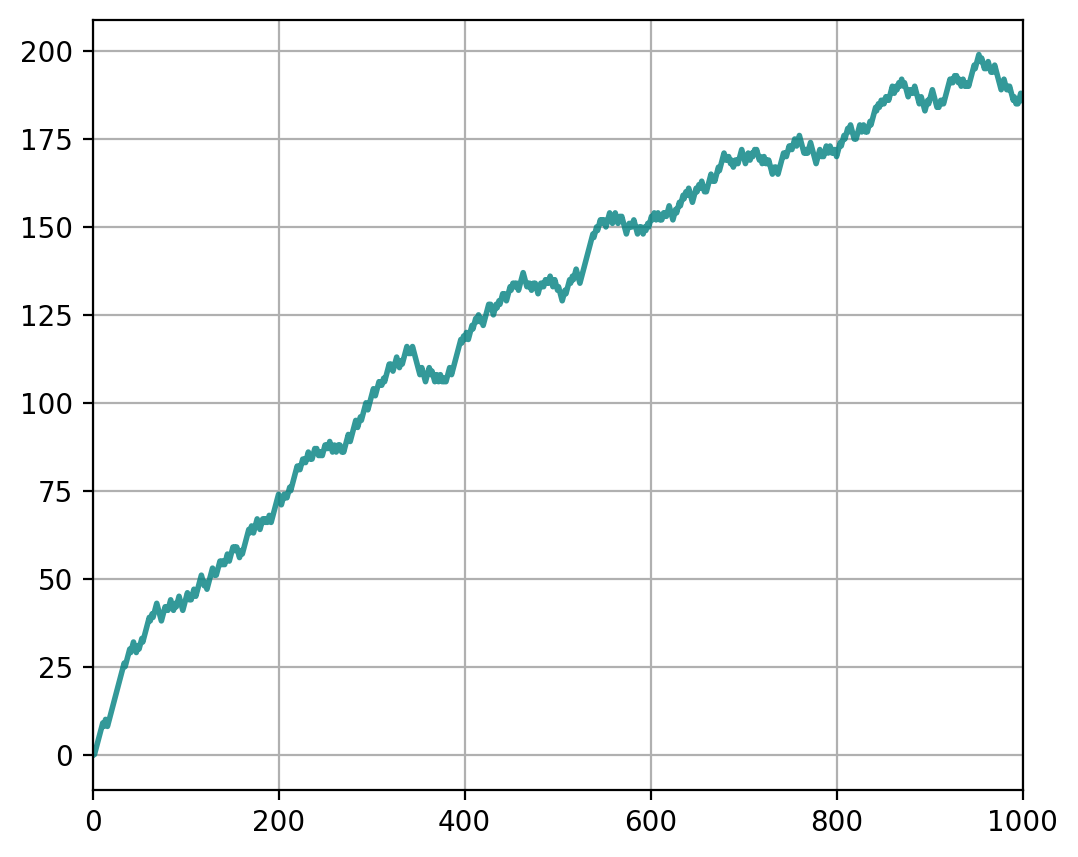}\qquad 
    \includegraphics[width=7cm,height=5.5cm]{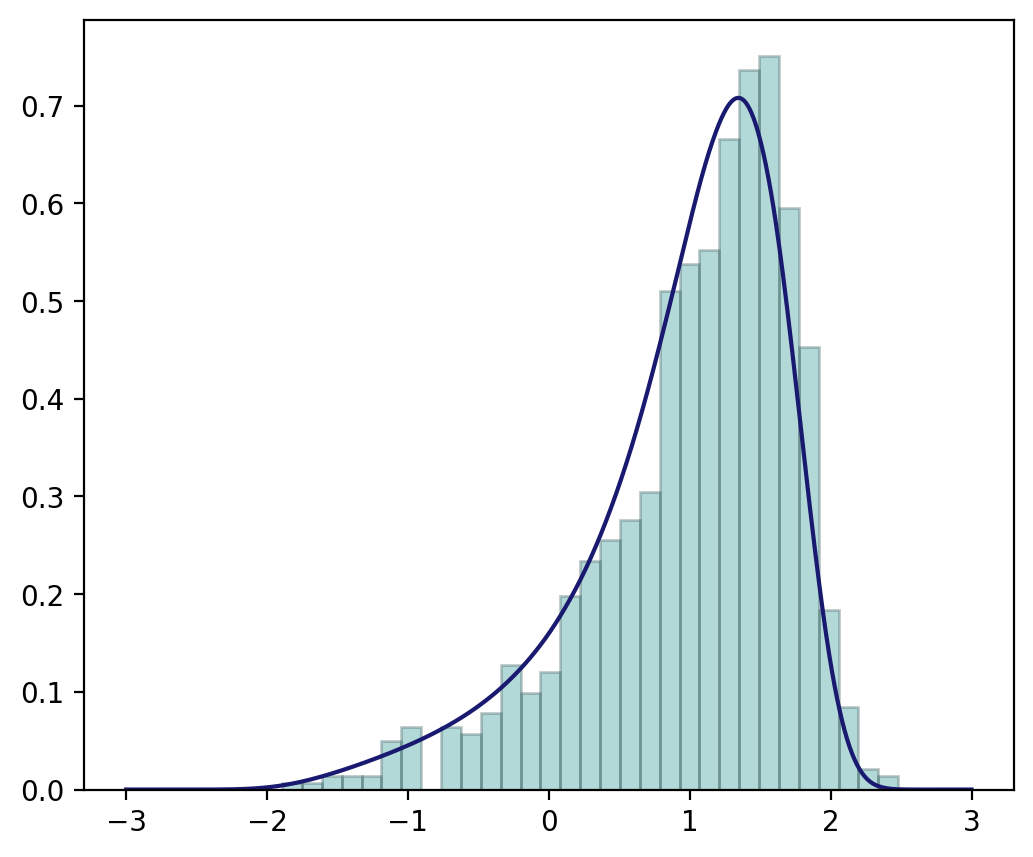}
    \caption{Left display: a trajectory of the elephant random walk in the superdiffusive regime until time $n=1000$, with $p=0.92$ and $q=1$. Right display: a simulation of the density of $L_1$ for the same parameters.} 
    \label{fig:simulations_ERW}
\end{figure}

The ERW process $(S_n)_{n\geq 0}$ is defined as follows. We denote its successive steps by $(X_n)_{n\geq 0}$.
The elephant starts at the origin at time zero: $S_0 = 0$. For the first step $X_1$, the elephant moves one step to the right with probability $q$ or one step to the left with probability $1-q$, for some $q$ in $[0,1]$. The next steps are performed by uniformly choosing an integer $k$ from the previous times. Then the elephant moves exactly in the same direction as at time $k$ with probability $p\in[0,1]$, or in the opposite direction with probability $1-p$.
In other words, defining for all $n \geq 1$, 
\begin{equation*}
   X_{n+1} = \left \{ \begin{array}{lll}
    +X_{\mathcal U(n)} &\text{with probability  $p$,} \\[1em]
    -X_{\mathcal U(n)} &\text{with probability $1-p$,}
   \end{array} \right.
\end{equation*}
with $\mathcal U(n)$ an independent uniform random variable on $\BRA{1,\ldots,n}$, the position of the ERW at time $n+1$ is 
\begin{equation*}
   S_{n+1}=S_{n}+X_{n+1}.
\end{equation*}
The probability $q$ is called the first step parameter and $p$ the memory parameter of the ERW. An ERW trajectory is sampled in Figure~\ref{fig:simulations_ERW}.

A wide range of literature is now available on the ERW and its extensions, see for instance
\cite{BaurBertoin2016,Coletti2017,Cressoni2013,Cressoni2007,daSilva2013,Kursten2016}. 
The behaviour of the process, in particular its dependency on the value of $p$ with respect to the critical value $3/4$, is now well understood. In the diffusive regime $p< 3/4$ and the critical regime $p=3/4$, a strong law of large numbers and a central limit theorem for the position, properly normalized, were established, see \cite{BaurBertoin2016,Coletti2017,ColettiN2017,Schutz2004} and the more recent contributions \cite{BercuHG2019,Coletti2019,Fan2020,Gonzales2020,Vazquez2019}. The main change between the two regimes is the rate of the associated convergences.

Introduce 
    $a:=2p-1$. In this paper, we will focus on the superdiffusive regime $a>1/2$ (i.e., $p>3/4$), which is more intriguing. 
It has been established that
\begin{equation*}
    \lim_{n \rightarrow \infty} \frac{S_n}{n^{a}}=L_q \quad \text{a.s.}
\end{equation*}
where $L_q$ is a non-degenerate, non-Gaussian random variable, see \cite{BaurBertoin2016,Bercu2018,Coletti2017}. Moreover, the fluctuations of the ERW around its limit $L_q$
are Gaussian by \cite{Kubota2019}.
These results were established thanks to a martingale approach.

While the limit variable $L_q$ has remained mysterious for some time, the recent paper \cite{GLR} established several properties of it, such as the existence of a density, a recursive formula for the moments, the finiteness of the moment-generating function and the moment problem. This was achieved by reformulating the ERW in terms of urn processes together with fixed point equations.

In the present work, we pursue this line of research and obtain fine properties of the distribution, via other methods. We discard the case $a=1$ which is of no interest, since then $L_q\overset{a.s.}{=}X_1$ is the first step of the elephant random walk. We focus the study on the distribution of $L_1$, the limit of the superdiffusive ERW random variable with parameter $a \in (1/2,1)$ and first step parameter $q=1$. As observed in \cite[Rem.~ 2.3]{GLR}, there is a simple relation between the distribution of $L_1$ and the general asymptotic variable $L_q$:
\begin{equation}
\label{Lq1}
L_q\,\overset{d}{=}\,(2\xi_q-1)L_1,    
\end{equation}
where $\xi_q$ is a Bernoulli variable with parameter $q$, independent of $L_1$. Our first result is the following. 

\begin{Theorem}
\label{Uni}
The random variable $L_1$ is unimodal.
\end{Theorem}

If we set $\varphi$ for the density function of $L_1$ on $\rl,$ which is known to exist thanks to \cite[Thm~1.3]{GLR}, let us recall that the unimodality of $L_1$ means that this function has a unique local maximum. Our method to obtain this important property is to show first the discrete unimodality of $S_n$ itself on $\dZ$ for all $n,$ and then take the renormalized limit defining $L_1$. This discrete unimodality is obtained by induction, using the (non-homogeneous) Markovian character of the ERW. A classical refinement of unimodality for sequences or functions is the log-concavity, which is also known as a strong unimodality property since it preserves unimodality under additive convolution. For more details on this topic, see for example the survey papers \cite{Saumard, Stan}. In the present paper, we are able to characterize the log-concavity of the ERW for all $n$ in terms of a certain threshold parameter $a_0$ -- see Proposition~\ref{LC}, which implies that the function $\varphi$ is log-concave on $\rl$ for all $a\in (1/2,a_0].$ We also believe that the ERW becomes log-concave for all $a\in (1/2,1)$ and $n$ large enough, but this non-trivial problem, which involves complicated polynomials, still eludes us.

\medskip

\begin{figure}
    \centering
    \includegraphics[height=4cm]{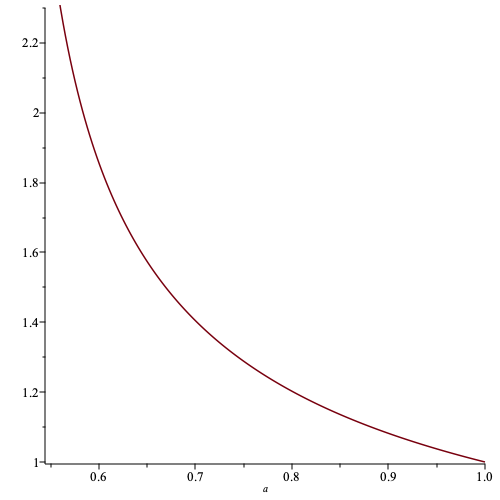}\qquad
    \includegraphics[height=4cm]{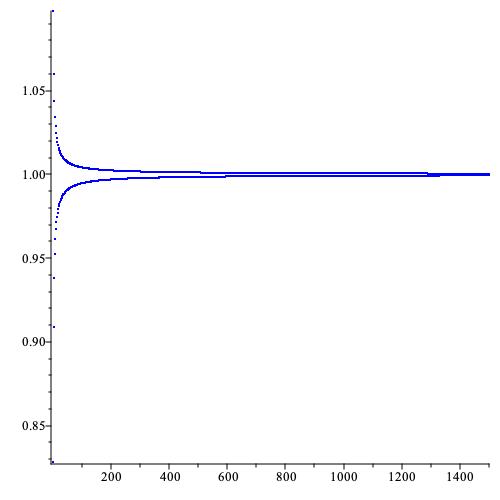}
    \qquad
    \includegraphics[height=4cm]{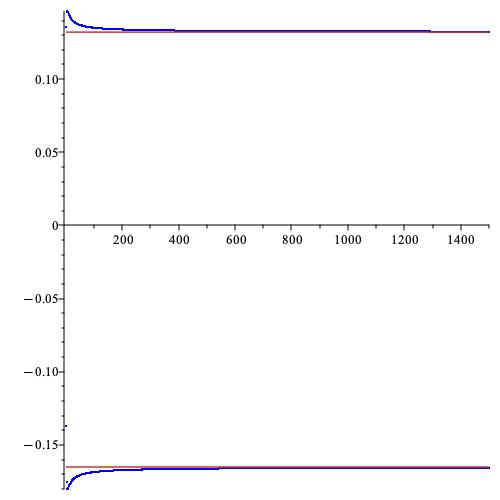}
    \caption{Left display: graph of the exponential growth $\rho_a$ on $(\frac{1}{2},1)$. The function $\rho_a$ is decreasing on $(\frac{1}{2},1)$, as will be proved in Proposition~\ref{ExpGro}. Middle display: the ratio $m_n/(\frac{2a}{a+1}\rho_a^n)$ converges to $1$ as $n\to+\infty$, according to Theorem~\ref{Mom} (computations made for $a=\frac{2}{3}$, using the exact recurrence relation \eqref{eq:recursion_moments}). Right display: the ratio $n^{\frac{2a}{a+1}}(m_n-\frac{2a}{a+1}\rho_a^n)/\rho_a^n$ converges to a constant depending on the parity of $n$, see \eqref{Asympto}, once again for $a=\frac{2}{3}$.}
    \label{fig:plot_rho}
\end{figure}

Our next results give the precise asymptotic behaviour of $\varphi$ at $\pm\infty$. These estimates involve a curious parameter $\rho_a$, which has not appeared yet in the ERW literature: 
\begin{equation}
    \label{eq:def_rho}
    \rho_a\, :=\, \lpa \frac{1}{\sqrt{\pi}} \, \Gamma\lpa \frac{1}{2} + \frac{1}{2a}\rpa \Gamma \lpa 1 - \frac{1}{2a}\rpa \rpa^a\, =\, \lpa\frac{{\rm B} \bigl(\frac{1}{2}+\frac{1}{2a},1-\frac{1}{2a}\bigl)}{2}\rpa^a.
\end{equation}
As proved in Proposition~\ref{ExpGro}, the function $a\mapsto \rho_a$ is decreasing on $(1/2,1)$ from $+\infty$ to $1$, see also Figure~\ref{fig:plot_rho}.

\begin{Theorem}
\label{Positano}
One has
\begin{equation*}
    \varphi(x)\, \sim\, c_a \, x^{\frac{2a-1}{2(1-a)}}\, e^{-(1-a)  \lpa \frac{a^a x}{\rho_a}\rpa^{\frac{1}{1-a}}}
\end{equation*}
as $x\to\infty$, with $c_a$ an explicit constant given in \eqref{eq:value_c_a}.
\end{Theorem}
For $a\in(1/2,1)$, we have $\frac{1}{1-a}>2$ and we hence recover the fact, first pointed out by Bercu in \cite[Rem.~3.2]{Bercu2004}, then refined in \cite[Cor.~2.12]{GLR}, that the distribution of $L_1$ is sub-Gaussian.

\begin{Theorem}
\label{Negus}
One has
\begin{equation*}
   \varphi(-x)\, \sim\, {\widehat c}_a \, x^{\frac{2a^2 -3a-1}{2(1-a^2)}}\, e^{-(1-a)  \lpa \frac{a^a x}{\rho_a}\rpa^{\frac{1}{1-a}}}
\end{equation*}
as $x\to\infty$, with $\widehat c_a$ an explicit constant given in \eqref{eq:value_c_a_hat}.
\end{Theorem}

Interestingly, the above estimates imply an asymmetry to the right of the random variable $L_1,$ at the level of the density function.
We notice that for $a\in(1/2,1)$ (actually for $a\in(0,1)$), $\frac{2a-1}{2(1-a)}>\frac{2a^2-3a-1}{2(1-a^2)}$. Due to the construction of $L_1$, it is not surprising to have a heavier tail at $+\infty$ compared to $-\infty$. In fact, $L_1$ is the asymptotic variable when the process starts from a step $X_1=+1$ a.s., with a memory parameter $p>3/4$. The memory parameter being large, the random walk tends to repeat its previous steps, with a first step to the right. The elephant is therefore more likely to move to the right. Observe that this asymmetry is peculiar to the case $q =1.$ For $q \in (0,1)$ indeed, the random variable $L_q$ has density $\varphi_q(x) = q\varphi(x) + (1-q)\varphi(-x)$ and the following corollary, which is an immediate consequence of Theorems~\ref{Positano} and \ref{Negus}, shows that this density has comparable tails at $\pm\infty.$

\begin{Cor}
\label{AsmpLq}
With the above notation, for $q\in(0,1)$, one has
\begin{equation*}
  \varphi_q(x)\, \sim\, q\, c_a \, x^{\frac{2a-1}{2(1-a)}}\, e^{-(1-a)  \lpa \frac{a^a x}{\rho_a}\rpa^{\frac{1}{1-a}}}\qquad \mbox{and}\qquad   \varphi_q(-x)\, \sim\, (1-q)\, c_a \, x^{\frac{2a-1}{2(1-a)}}\, e^{-(1-a)  \lpa \frac{a^a x}{\rho_a}\rpa^{\frac{1}{1-a}}}
\end{equation*}
as $x\to\infty$.
\end{Cor}

In the case $q= 1$, the difference between the two tails of the densities is specified in the next result, which is also a direct consequence of Theorems~\ref{Positano} and \ref{Negus}. 
 
\begin{Cor}
\label{Asym}
With the above notation, one has
\begin{equation*}
   \frac{\varphi(x)}{\varphi(-x)}\, \sim\, \frac{4\,\Gamma\lpa\frac{1-a}{1+a}\rpa}{(a+1)^{\frac{2a}{1+a}}} \lpa \frac{a}{\rho_a^{2-1/a}}\rpa^{\frac{2a}{1-a^2}}\, x^{\frac{2a}{1-a^2}}\qquad \mbox{as $x\to\infty.$}
\end{equation*}
\end{Cor}

Our proofs of Theorems \ref{Positano} and \ref{Negus} are based on a combination of exact computations using special functions (such as hypergeometric functions,  incomplete beta integrals and Mittag-Leffler functions) together with singularity analysis of generating functions (such as strong Tauberian theorems). Although the techniques for proving the two tails of the density are similar, the asymptotics on the half-negative line is significantly trickier. Let us present three main aspects of our methods.

\medskip

The first step is to start from the moment calculation made in \cite[Thm~1.4]{GLR}
\begin{equation}
\label{eq:link_moments_m_n}
   \dE[L_1^n]\,=\,\frac{n!\, m_n}{\Ga(1+an)},\qquad \mbox{for all $n\geq 1$,}
\end{equation}
and the non-linear recurrence relation for the sequence $\{m_n\}$ to be recalled in \eqref{eq:recursion_moments}. We will deduce\footnote{In general, if a power series \eqref{eq:def_M} satisfies an algebraic differential equation, then its coefficients $\{m_n\}$ should satisfy a recurrence relation. On the other hand, the existence of a recurrence relation at the level of the coefficients does not necessarily imply that the associated power series satisfies a differential equation. In other words, the existence of a non-trivial differential equation for our power series $M(x)$ in \eqref{eq:def_M} is a result in itself.} an algebraic differential equation for the Stieltjes-like generating function
\begin{equation}
\label{eq:def_M}
   M(x)\,=\,\sum_{n\geq 0}m_n x^n
   \,=\, 1\, +\, x \, +\, \frac{a}{2a-1}x^2\,+\,\frac{a+1}{2(2a-1)}x^3\,+\cdots,
\end{equation}
see \eqref{eq:ODE_delay_M}. Surprisingly, despite its apparent complexity, the differential equation \eqref{eq:ODE_delay_M} admits a closed-form solution, which can be expressed in terms of a certain (incomplete beta)\ integrals, or equivalently in terms of some hypergeometric functions; see Proposition~\ref{prop:exact}. In a sense, the explicit expression for $M(x)$ is an integrability property of the ERW.

In the second step, we prove that the series \eqref{eq:def_M} has a first positive singularity at the point $\frac{1}{\rho_a}$, with $\rho_a$ introduced in \eqref{eq:def_rho}, and we obtain a series expansion at that point with the help of the exact expression given in the first step. Using the singularity analysis on generating functions developed in \cite{FO}, we can deduce an asymptotic expansion of the moments, viewed as the coefficients of the power series \eqref{eq:def_M}. See Theorem~\ref{Mom} for the precise statement.

The third step aims at proving Theorems~\ref{Positano} and \ref{Negus} on the asymptotic behavior of the density $\varphi$ at $\pm\infty$. The general idea is to apply the strong Tauberian theorems of \cite{FY}, where the fine properties of the exponential generating function
\begin{equation}
\label{eq:def_MGF}
   \Psi(r)\, =\, \esp[e^{rL_1}]\, =\, \sum_{n\geq 0} \lpa\frac{\esp[L_1^n]}{n!}\rpa r^n
\end{equation}
can be transferred into asymptotic estimates of the density at $\pm\infty$. This is the most delicate part of the paper, which intrinsically requires an expansion of the integer moments up to the fifth order. The latter implies an expansion of $\Psi(r)$ in terms of generalized Mittag-Leffler functions and we can then apply the global estimates of Wright -- see \cite{W40,GG} -- on such functions in order to derive a sufficient control on $\Psi$ and its derivatives and to perform accurately the inverse Laplace approximation method of \cite{FY}. For the reader's comfort, we also offer as an intermediate result in Proposition \ref{Log} a proof of the asymptotic behaviour at both infinities under the logarithmic scale, which is a simpler consequence of Kasahara's exponential Tauberian theorem. 

\section{Proof of Theorem \ref{Uni}}

By definition, we have 
\begin{equation*}
   n^{-a} S_n \, \claw\, L_1
\end{equation*}
as $n\to\infty.$ We first show that the random variable $S_n,$ taking values in $\{2k -n, \, k = 1,\ldots, n\},$ is unimodal for every $n\geq 1$ and all $p\in [0,1]$ (which covers cases where the ERW is not necessarily superdiffusive). 

\medskip

Setting $P(n,k) = \pb[S_n = 2k-n]$ and $Q(n,k) = (n-1)! P(n,k),$ it follows from the definition of the ERW that
$\{Q(n,k), \, 1\leq k \leq n,\, n\geq 1\}$ is a triangular array defined recursively by $Q(1,1) = 1$ and
\begin{equation}
\label{Recursive}
Q(n+1,k) \, =\, (np -ak) Q(n,k) \, +\, (qn + a(k-1)) Q(n,k-1)
\end{equation}
for all $n\geq 1$ and $1\leq k\leq n+1,$ where we have set $q= 1-p$ and, here and throughout, we make the convention $Q(n,0) = Q(n,n+1) = 0.$ 
Let us briefly prove \eqref{Recursive}. The Markov property yields 
\begin{equation*}
    P(n+1,k) = \mathbb P[S_{n+1}=2k-n-1\vert S_n=2k-n]P(n,k)+\mathbb P[S_{n+1}=2k-n-1\vert S_n=2k-n-2]P(n,k-1).
\end{equation*}
To conclude to \eqref{Recursive}, we multiply the above equation by $n!$, use the notation $Q(n,k)$ and the explicit transition probabilities of the ERW, namely 
\begin{equation*}
    \left\{\begin{array}{lcl}
    \displaystyle\mathbb P[S_{n+1}=2k-n-1\vert S_n=2k-n] &=&\displaystyle\frac{1}{2}\left(1-a\frac{2k-n}{n}\right),\smallskip\\
    \displaystyle\mathbb P[S_{n+1}=2k-n-1\vert S_n=2k-n-2] &=&\displaystyle\frac{1}{2}\left(1+a\frac{2k-n-2}{n}\right).
    \end{array}\right.
\end{equation*}

We need to show that for every $n\geq 1$ the finite sequence $\{Q(n,k), \, 1\leq k \leq n\}$ is unimodal, that is -- see e.g.\ Definition~4.1 in \cite{DJD} -- there exists a maximal index $k_n \in \{1,\ldots, n\}$ such that
\begin{equation*}
   \lacc\begin{array}{ll}
Q(n,k) \geq Q(n, k-1) & \mbox{for all $k\leq k_n$}, \\
Q(n,k) \leq Q(n, k-1) & \mbox{for all $k > k_n.$}
\end{array}
\right.
\end{equation*}
To do so, we use an induction on $n$ starting from the unimodal sequence $\{Q(1,1) = 1\}.$ Let $n\geq 1$ and suppose that the finite sequence $\{Q(n,k), \, 1\leq k\leq n\}$ is unimodal. Decomposing
\begin{multline*}
    Q(n+1,k+1)\, -\, Q(n+1,k)=\\
    (np-a(k+1))\bigl(Q(n,k+1)-Q(n,k)\bigr) \, +\, (nq+a(k-1)) \bigl(Q(n,k) -Q(n,k-1)\bigr)
\end{multline*}
implies that $Q(n+1,k+1) \geq Q(n+1,k)$ for all $k = 1,\ldots, k_n-1$ by the induction hypothesis, since then one has $np - a(k+1) \geq np - ak_n \geq n(p-a) = n(1-p)\geq 0.$ In particular, if $k_n = n,$ the sequence $\{Q(n+1,k), \, 1\leq k\leq n+1\}$ is unimodal with maximal index $k_{n+1} = n$ or $k_{n+1} = n+1.$ In the case $k_n < n,$ one necessarily has $Q(n,n) \leq Q(n, n-1)$ by the maximality of $k_n,$ and we first show that this implies
\begin{equation}
\label{pmax}
p \,\leq\, \frac{n+1}{n+2}\cdot
\end{equation}
Indeed, if the contrary held, then we would also have $(k+2)p > k+1\Leftrightarrow a(k+1) > kp$ for all $k = 1,\ldots, n.$ Since on the other hand one has
\begin{equation*}
   Q(k+1,k+1)  -  Q(k+1,k) \, = \, (a(k+1) - kp) Q(k,k) \, +\, (kq + a(k-1)) (Q(k,k)  - Q(k,k-1))
\end{equation*}
for all $k=1,\ldots,n,$ an induction starting from $Q(1,0) = 0$ and $Q(1,1) = 1$ readily implies then $Q(k+1,k+1) > Q(k+1,k)$ for all $k = 1, \ldots, n,$ a contradiction in the case $k = n-1.$ 

Therefore, since \eqref{pmax} implies $np \geq a(k+1)$ for all $k = 1, \ldots, n,$ we also have $Q(n+1,k+1) \leq Q(n+1,k)$ for all $k = k_n+1, \ldots, n.$ Putting everything together, we deduce that the sequence $\{Q(n+1,k), \, 1\leq k\leq n+1\}$ is unimodal as required, with maximal index $k_{n+1} = k_n$ if $Q(n+1,k_n+1) < Q(n+1,k_n)$ and $k_{n+1} > k_n$ if $Q(n+1, k_n +1) \geq Q(n+1,k_n).$  

\medskip

In order to complete the proof of the unimodality of $L_1$, we focus on the superdiffusive case $a\in (1/2,1)$ and consider the step function $f_n$ defined by 
\begin{equation*}
   f_n(x)\, = \, \frac{n^{a} P(n,k)}{2}
\end{equation*}
if $x\,\in\, (n^{-a}(2k-n-1), n^{-a}(2k-n+1)]$ for some $k\in\{1,\ldots, n\}$ and $f_n(x) = 0$ otherwise. It is clear from the above unimodality of $S_n$  that $f_n(x)$ is a unimodal density function on $\rl.$ Moreover, it is easy to see by construction that the corresponding distribution function $F_n (x) = \int_{-\infty}^x f_n(y)\, dy$ is such that
\begin{equation*}
   \sup_{x\in\rl} \lva F_n(x)\, -\, \pb[n^{-a} S_n \leq x]\rva \, \leq\, \max\{P(n,k), \, k=1,\ldots, n \}.
\end{equation*}
Setting $m_n = \max\{Q(n,k), \, k=1, \ldots, n-1\}$ and
\begin{equation*}
   M_n\, =\, \max\{Q(n,k), \, k=1, \ldots, n\}\, =\, \max\{(n-1)! p^{n-1}, m_n\}
\end{equation*}
for all $n\geq 2,$ we see by definition that $m_{n+1} \leq (n-a)M_n$ and
$Q(n+1,n+1) = pn Q(n,n) \leq (n-a) M_n$ for $n$ large enough. This shows that there exists $n_0\geq 1$ such that $M_{n+1} \leq (n-a) M_n$ for all $n\geq n_0,$ whence
\begin{equation*}
   \max\{P(n,k), \, k=1,\ldots, n \} \, =\, \frac{M_n}{(n-1)!}\, \leq\, \frac{M_{n_0}}{(1-a)_{n_0-1}}\lpa\frac{(1-a)_{n-1}}{(n-1)!}\rpa\,\to\, 0\quad\mbox{as $n\to \infty.$}
\end{equation*}
Putting everything together shows that $F_n(x)\to \pb[L_1\leq x]$ for all $x\in\rl$ and, by stability of unimodality under weak convergence -- see e.g.\ Theorem 1.1 in \cite{DJD}, that $L_1$ is unimodal. 
\qed

\begin{Rem}{\em (a) If $\varphi'(0) > 0,$ then the modes of $L_1$ are necessarily  positive. Observe that for the Mittag-Leffler random variable with parameter $a > 1/2,$ which is an approximation of $L_1$ from the point of view of integer moments and is also a strictly unimodal random variable, the mode is positive (whereas it is zero if $a\leq 1/2$). The known lower bound $M > \esp[L_1] - \sqrt{3}{\rm Var}[L_1]$ where $M$ is any mode of $L_1$ -- see e.g.\ Lemma~1.7 in \cite{DJD} implies that these modes are indeed positive for $a$ close enough to 1.

(b) Simulations show that in the case $q\in(0,1),$ the random variable $L_q$ is not unimodal in general. By \eqref{Lq1} and Theorem \ref{Uni} the density of $L_q$ is the superposition of two unimodal functions and we conjecture that for every $a\in (1/2,1)$ there exists $q(a)\in[0,1/2]$ such that $L_q$ is bimodal if $\vert q-1/2\vert < q(a)$ and unimodal if $\vert q-1/2\vert \ge q(a).$ This problem is believed to be challenging. }
\end{Rem}

We now show a refinement of the unimodality property for the ERW for certain values of the parameter $a.$ Recall that a sequence $\{u_k, \,k = 1,\ldots, n\}$ of positive numbers is log-concave if 
$$u_k^2 \,\geq\, u_{k-1}u_{k+1}$$ 
for all $k = 1, \ldots, n,$ where we have set $u_0 = u_{n+1} = 0.$  

\begin{Prop}
\label{LC}
The sequence $\{\pb[S_n = 2k-n], \, k=1,\ldots, n\}$ is log-concave for every $n\geq 1$ if and only if $a\in[-1,a_0]$ where $a_0 = 0.61803...$ is the unique positive root of $a^3 + 4a^2+2a = 3.$   
\end{Prop}

\proof With the above notation, we have to show that the sequence $\{Q(n,k), \, k=1,\ldots, n\}$ is log-concave for all $n\geq 1$ if and only if $a\in [-1,a_0].$ The sequences $\{Q(1,1)\}$ and $\{Q(2,1), Q(2,2)\}$ are always log-concave. We compute $Q(3,1)= (1-a)/2, Q(3,2) = (1-a)(2+a)/2, Q(3,3) = (1+a)^2/2$ and see that the log-concavity of $\{Q(3,1), Q(3,2), Q(3,3)\}$ amounts to
\begin{equation*}
   (1-a)(2+a)^2 \,\geq\, (1+a)^2\,\Longleftrightarrow\, 3\,-\, 2a\, -\, 4a^2\, -\, a^3 \,\geq \, 0\, \Longleftrightarrow\, a \in [-1,a_0],
\end{equation*}
which shows the only if part. We next suppose $a\in [-1,a_0]$ and show that $\{Q(n,k), \, k=1,\ldots, n\}$ is log-concave for all $n\geq 1$ by an induction on $n.$ From the above discussion, the property is true for $n=1,2,3$ and we will show later that it is also true for $n = 4.$ Assuming that $\{Q(n,k), \, k=1,\ldots, n\}$ is log-concave for some $n\geq 4,$ we first show 
\begin{equation*}
   Q(n+1,n)^2\, \geq\, Q(n+1,n+1)Q(n+1,n-1).
\end{equation*}
Setting $A_n = Q(n+1,n)^2\, -\, Q(n+1,n+1)Q(n+1,n-1)$ and using \eqref{Recursive} we compute, after some simplifications, 
\begin{eqnarray*}
4 A_n & = & n^2(1-a)^2Q(n,n)^2\,  + \, n(n(1-a^2) -6a +2a^2)Q(n,n)Q(n,n-1)\,+\, 4a^2Q(n,n-1)^2\\
& & \,+\; n(1+a)(n(1+a) -4a))(Q(n,n-1)^2\, -\, Q(n,n)Q(n,n-2))\\
& \geq & n^2(1-a)^2Q(n,n)^2\,  + \, n(4 -6a - 2a^2) Q(n,n)Q(n,n-1)\,+\; 4a^2Q(n,n-1)^2\\
& \geq & n^2(1-a)^2Q(n,n)^2\,  - \, 4na(1-a)Q(n,n)Q(n,n-1)\,+\; 4a^2Q(n,n-1)^2\\
& = & \lpa n(1-a) Q(n,n) \, -\, 2a Q(n,n-1)\rpa^2\, \geq \, 0,
\end{eqnarray*}
where in the first inequality we have used the induction hypothesis and $n\geq 4,$ and in the second inequality we have used $a\leq a_0 \leq 2/3$ which implies $4 -6a -2a^2 + 4a(1-a) = 2(1+a)(2-3a)\geq 0.$ We next prove
\begin{equation*}
   Q(n+1,k)^2\, \geq\, Q(n+1,k+1)Q(n+1,k-1)
\end{equation*}
for all $k=2,\ldots,n-1.$ Setting $B_{k,n} = Q(n+1,k)^2 - Q(n+1,k+1)Q(n+1,k-1)$ and using \eqref{Recursive} again we compute
\begin{eqnarray*}
B_{k,n} & = & (np -a(k+1))(np -a(k-1)) B_{k,n-1} \, +\, (nq+ak)(nq +a(k-2) B_{k-1,n-1}\\
& & \, +\; a^2\lpa Q(n,k) - Q(n,k-1)\rpa^2 \, -\, (np-a(k+1))(nq+a(k-2)) Q(n,k-2)Q(n,k+1)\\ 
& & \, +\; \lpa 2a^2 + 2(np-ak)(nq +a(k-1)) - (nq+ak)(np - a(k-1))\rpa Q(n,k-1)Q(n,k)\\
& \geq & (np-a(k+1))(nq+a(k-2))\lpa Q(n,k-1)Q(n,k) - Q(n,k-2)Q(n,k+1)\rpa\, \geq \, 0,
\end{eqnarray*}
where the second equality follows from the induction hypothesis. To complete the proof, it remains to show that $\{Q(n,k), \, k=1,\ldots, n\}$ is log-concave for all $a\in [-1,a_0]$ and $n=4.$ The above discussion and the case $n= 3$ already implies $Q(4,2)^2 \geq Q(4,1)Q(4,3)$ and we hence just need to check that $A_3 \geq 0,$ with the above notation, which amounts after some simplifications to
\begin{equation*}
   P(a) \, =\, 54\, +\, 36a \, - 83a^2\, -\, 123 a^3\, -\, 63a^4 \, -\, 13 a^5 \, \geq \, 0.
\end{equation*}
The polynomial $P(a)$ is easily seen to decrease on $[1/2,1]$ with $P(1/2) > 0$ and $P(1) < 0$ and has hence a unique root $a_1 = 0.63606.. > a_0,$ so that $P(a) \geq 0$ for all $a\in [1/2,a_0].$ Moreover, we have 
\begin{eqnarray*}
4 A_3 & \geq & 9(1-a)^2Q(3,3)^2\,  + \, 3(3 -6a - a^2) Q(3,3)Q(3,2)\,+\; 4a^2Q(3,2)^2\\
& \geq & \lpa 3(1-a) Q(3,3) \, -\, 2a Q(3,2)\rpa^2\, \geq \, 0
\end{eqnarray*}
as well for all $a\in[-1,1/2],$ since then $3-6a-a^2 + 4a(1-a) = (1+a)(3-5a)\geq 0.$ 
\qed

\begin{Rem}
\label{LC1}
{\em One can show that there exists an increasing sequence $\{a_q, q\geq 0\}$ such that
\begin{equation*}
   \{Q(n,k), \, k=1,\ldots,n\}\; \mbox{is log-concave for all $n\geq q+3$}\;\Longleftrightarrow\; a\,\leq\, a_q
\end{equation*}
for all $q\geq 0,$ and that $x_q$ is the unique root on $(1/2,1)$ of a certain real polynomial $P_q$ whose degree is smaller than $2q+3.$ One has

\medskip

$$\begin{array}{l}
P_0 \, =\, 3\, -\, 2X\, -4X^2\, - \, X^3\\
\\
P_1 \, =\, 54\, +\, 36X \, - \, 83 X^2\, -\, 123 X^3\, -\, 63X^4 \, -\, 13 X^5
\\
\\
P_2 \, =\, 360\, +\, 834X \, +\, 247 X^2\, -\, 1259 X^3\, -\, 2009 X^4 \, -\, 1423 X^5\, -\, 514 X^6 \, -\, 76X^7
\end{array}$$
and 
$$P_3  = 11250 + 8325X  - 13781 X^2- 24282 X^3 - 13396 X^4 - 1024 X^5 + 2518X^6  + 1361X^7  + 229X^8$$
with $a_0 \simeq 0.61803, a_1 \simeq 0.63606, a_2\simeq 0.67060$ and $a_3\simeq 0.68408.$ We conjecture that 
\begin{equation*}
   \lim_{q\to\infty} a_q\, = \,1.
\end{equation*}
However, this problem seems challenging because of the very complicated character of the underlying polynomials.}
\end{Rem}

The following is an interesting consequence of Proposition \ref{LC}.

\begin{Cor}
\label{LC2}
The density of $L_1$ is log-concave on $\rl$ for all $a\in(1/2,a_0].$ 
\end{Cor}

\proof The argument is analogous to that of the end of the proof of Theorem \ref{Uni} except that we consider here the piecewise affine density $g_n$ defined for all $n\ge 2$ by

\begin{equation*}
   g_n(x)\, = \, \frac{n^a}{2- (P(1,n) + P(n,n))}\lpa P(k,n) + \frac{n^a(P(k+1,n) - P(k,n))(x + n - 2k)}{2}\rpa
\end{equation*}
if $x\,\in\, (n^{-a}(2k-n), n^{-a}(2k-n+2)]$ for some $k\in\{1,\ldots, n-1\}$ and $g_n(x) = 0$ otherwise. The log-concavity of the sequence $\{P(n,k), \, k=1, \ldots, n\}$ readily implies the log-concavity of the function $g_n$ on $\rl$ for all $a\in (1/2,a_0]$, and it follows easily from Scheffé's lemma and the end of the proof of Theorem \ref{Uni} that the real random variable with density $g_n$ converges weakly to $L_1$ as $n\to\infty$. Since log-concavity is a stable property under weak convergence -- see e.g.\ Theorem~2.10 in \cite{DJD}, this shows altogether that $L_1$ has a log-concave density for all $a\in (1/2,a_0]$.
\qed

\begin{Rem}
\label{LC3}
{\em A positive answer to the open problem stated in Remark \ref{LC1} would imply by the same approximation argument that $L_1$ has a log-concave density on $\rl$ for all $a\in (1/2,1).$ Observe that this property is in accordance with the superexponential behaviour of the density at $\pm\infty$ stated in Theorems  \ref{Positano} and \ref{Negus}.}
\end{Rem}

\section{Moment estimates}
In this section we establish sharp estimates for the positive integer moments of $L_1,$ which play a crucial role in the proof of Theorems \ref{Positano} and \ref{Negus}, and have an independent interest. The principal estimate is given in \eqref{Asympto1} below. However, the proof of Theorem \ref{Positano} requires an expansion up to the third order whereas that of Theorem \ref{Negus} needs an expansion up to the fifth order, which is given in \eqref{Asympto} with the notation
\begin{equation}
\label{eq:def_kappa_a_delta_a}
   \kappa_a\, =\, \frac{\rho_a^{\frac{2}{a+1}}}{4}\lpa \frac{a+1}{a}\rpa^{\frac{2a}{a+1}}\qquad\mbox{and}\qquad \delta_a \, =\, \frac{1-a}{1+a}\cdot
\end{equation}
Our main tool to get the moment estimates is the singularity analysis of \cite{FO}, where expansions are given in terms of Pochhammer symbols which we recall to be defined by $(x)_0 = 1$ and $(x)_n = x(x+1)\cdots (x+n-1)$, for all $n\geq 1$ and $x\in\rl.$ This can be easily transformed into an expansion in descending powers of $n$ thanks to the asymptotics 
\begin{equation}
\label{Posh}
\frac{(x)_n}{n!}\, =\, \frac{\Gamma(x+ n)}{\Gamma(x) \Gamma (1+n)}\, =\, \frac{n^{x-1}}{\Gamma(x)}\lpa 1\, +\, \frac{x(x-1)}{2n}\, +\, {\rm O} \lpa n^{-2}\rpa\rpa,
\end{equation}  
which is valid for all $x > 0$ -- see \cite{ET} for a more complete version. The asymptotic expansion with Pochhammer symbols will also be convenient along the proof of Theorems \ref{Positano} and \ref{Negus}, thanks to the direct connection with generalized Mittag-Leffler functions and their behaviour at infinity.
 
\begin{Theorem}
\label{Mom}
One has
\begin{equation}
\label{Asympto1}
\esp[L_1^n]\, \sim \, \frac{2a\, \rho_a^n \, n!}{(a+1)\,\Ga(1+an)}\qquad\mbox{as $n\to\infty.$}
\end{equation}
More precisely, one has the asymptotic expansion
\begin{multline}
\label{Asympto}
  \esp[L_1^n]\, = \, \frac{2a\, \rho_a^n \, n!}{(a+1)\,\Ga(1+an)}\lpa 1\, +\, \kappa_a\, \frac{(\delta_a)_n}{n!}\lpa (-1)^n\, + \, \frac{a-1}{3a+1}\rpa\right.\\ \left.\, +\, \frac{(2\delta_a)_n}{(n+1)!}\lpa (-1)^n \kappa_1^- + \kappa_1^+\rpa \, + \, \frac{(\delta_a)_n}{(n+1)!}\lpa (-1)^n \kappa_2^- + \kappa_2^+\rpa\, +\,{\rm O} \lpa n^{-2}\rpa\rpa 
\end{multline}
as $n\to\infty,$ where $\rho_a$ is introduced in \eqref{eq:def_rho}, $\kappa_a$ and $\delta_a$ in \eqref{eq:def_kappa_a_delta_a}, $\kappa_1^-, \kappa_1^+, \kappa_2^-$ and $\kappa_2^+$ are computable real constants depending on $a$.
\end{Theorem}

The proof of this theorem relies on the representation \eqref{eq:link_moments_m_n}, where the sequence $\{m_n\}$ is defined by $m_0= m_1 = 1$ and
\begin{equation}
\label{eq:recursion_moments}
m_n\,=\,\frac{1}{na-c_n}\sum_{i=1}^{n-1}c_im_im_{n-i},
\end{equation}
for every $n\geq 2,$ where $c_i=1$ for even $i$ and $c_i=a$ for odd $i$. By Lemma~2.10 and Corollary~2.14 in \cite{GLR}, for every $a\in\PAR{\frac{1}{2},1}$ there exists $C_a > 0$ such that $0 < m_n\leq C_a^n$ for all $n\geq 0,$ which shows that the generating function \eqref{eq:def_M}
has a positive radius of convergence. We will first give an exact expression of $M(x),$ showing  that this radius is actually $1/\rho_a$, with $\rho_a$ introduced in \eqref{eq:def_rho}. Set
\begin{equation}
\label{eq:def_F_cv}
   F(x)\, =\, \frac{1}{a}\int_x^\infty \frac{du}{u^{1+1/a}\sqrt{1+u^2}}
\end{equation}
and consider the function 
\begin{equation}\label{eq:def_G}
    G(x)\,=\, F^{-1}(x^{-1/a}- \rho_a^{1/a})
\end{equation} 
for every $x\in(0,1/\rho_a),$ where $F^{-1}$ denotes the compositional inverse of the decreasing function $F$ on $(0,\infty).$ 

\begin{Prop}
\label{prop:exact}
For every $x\in(0, 1/\rho_a)$, one has 
\begin{equation*}
\label{eq:expression_M_F}
    M(x)\, =\,  \left(\frac{G(x)}{x}\right)^{\frac{1}{a}}\left(G(x) + \sqrt{1+G(x)^2} \right),
\end{equation*}
with $G$ defined by \eqref{eq:def_G}.
\end{Prop}

The next paragraph is devoted to the proof of this proposition. Then, in Section~\ref{subsec:proof_as}, we will analyze precisely the behaviour of $M(x)$ as $x\to 1/\rho_a$ and prove Theorem \ref{Mom}.
 
\subsection{Proof of Proposition~\ref{prop:exact}} 
Multiplying the recurrence relation \eqref{eq:recursion_moments} by $(na-c_n) x^n$ and taking the generating functions yields first the following ODE for the series $M(x)$ in \eqref{eq:expression_M_F}: 
\begin{equation}
    \label{eq:ODE_delay_M}
    M(x)\, + \, axM'(x)\, - \, pM(x)^2 \, -\,  (1 - p)M(x)M(-x)\,=\, 0,
\end{equation}
which appears as a mix of an ordinary (algebraic) differential equation and a functional equation, because of the term $M(-x)$. Separating $M(x)=A(x)+B(x)$ with 
\begin{equation*}
   A(x)\, =\, \sum_{n\geq 0} m_{2n} x^{2n} \qquad\mbox{and} \qquad B(x)\, =\, \sum_{n\geq 0} m_{2n+1} x^{2n+1},
\end{equation*}
the above equation can be rewritten as
\begin{equation*}
   A\,+\,B\,+\,axA'\, +\, axB'\, -\, p(A^2+2AB+B^2)\, -\, (1-p)(A^2-B^2) \,= \, 0
\end{equation*}
and isolating the even and odd parts, we obtain the differential system
\begin{equation}
\label{eq:diff_system}
   \left\{\begin{array}{l}
    A\, +\, axA'\, -\, p(A^2+B^2)\, -\, (1-p)(A^2-B^2) \, =\, 0,\smallskip\\
    B\, +\, axB'\, -\, 2pAB\, =\, 0.
    \end{array}\right.
\end{equation}
Using the second equation we express $A$ in terms of $B$ and $B'$, and then plugging this identity in the first equation, we obtain after cleaning denominators and simplification an autonomous, classical ODE on $B(x)$, which reads 
\begin{equation}
\label{eq:diff_eq_B}
    a(a+1)x^2BB''-a(a + 2)x^2 B'^2 + x((a+1)^2 - 2)BB'- (a+1)^2 B^4     + B^2\, = \,0.
\end{equation}
The following lemma shows that $B(x)$ is also the solution to an implicit equation, involving the function $F$ in \eqref{eq:def_F_cv}. It is worth noticing that Lemma~\ref{prop:explicit_implicit_B} is the first place where the non-trivial quantity $\rho_a$ appears.
\begin{Lemma}
\label{prop:explicit_implicit_B}
For every $x\in (0,1/\rho_a)$, the function $B(x)$ is the unique solution of
\begin{equation}
\label{eq:implicit_B}
 x^{1/a} F\lpa (x^{1/a}B(x))^{a/(a+1)}\rpa\, + \,\PAR{\rho_a x}^{1/a}\, =\, 1.
\end{equation}
In particular, one has $B(x)\to\infty$ as $x\to 1/\rho_a$ and the convergence radius of the series defining $B(x)$ is $1/\rho_a.$ 
\end{Lemma}

\proof We first remark that the solution of \eqref{eq:implicit_B} exists and is uniquely defined as
\begin{equation}
\label{eq:expression_B_M_F}
    B(x) \, =\,  x^{-1/a}\left(F^{-1}(x^{-1/a}- \rho_a^{1/a})\right)^{\frac{a+1}{a}},
\end{equation}
which is a smooth function on $(0,1/\rho_a).$ Moreover, since  
\begin{equation*}
   x^{1/a} F\lpa (x^{1/a}B(x))^{a/(a+1)}\rpa\, = \, \frac{x^{1/a}}{a}\int_{(x^{1/a}B(x))^{a/(a+1)}}^\infty \frac{\dd u}{u^{1+1/a}}\PAR{\frac{1}{\sqrt{1+u^2}} -1}\,+\,\PAR{\frac{x}{B(x)}}^{1/(a+1)}
\end{equation*}
and since for any $a\in(1/2,1)$, the integral on the right-hand side is bounded as a function of $x$, we obtain $\lim_{x\to 0}\frac{x}{B(x)}= 1,$ whence $B$ is right-differentiable at zero with $B(0)=0$ and $B'(0)=1.$

We now prove that $B$ satisfies \eqref{eq:diff_eq_B}. Dividing \eqref{eq:implicit_B} by $x^{1/a}$ and differentiating with respect to $x$ we obtain, after a few elementary computations, 
    \begin{equation*}
        \frac{1}{a+1}\left(1+ax\frac{B'}{B}\right)=\left(\frac{B}{x}\right)^{\frac{1}{a+1}}\sqrt{1+x^{\frac{2}{a+1}}B^{\frac{2a}{a+1}}}.
    \end{equation*}
    Taking the square of the previous identity and simplifying, we get
    \begin{equation}
    \label{eq:ident_squar}
        \left(\frac{B}{x}\right)^{\frac{2}{a+1}} = -B^2+\frac{1}{(a+1)^2}\left(1+ax\frac{B'}{B}\right)^2.
    \end{equation}
    Differentiating \eqref{eq:ident_squar} with respect to $x$ and multiplying by $\frac{B}{x}$, we obtain
    \begin{equation}
    \label{eq:ident_squar_diff}
        \frac{2}{a+1}\left(\frac{B}{x}\right)'\left(\frac{B}{x}\right)^{\frac{2}{a+1}} =\frac{B}{x} \frac{\dd }{\dd x}\left(-B^2+\frac{1}{(a+1)^2}\left(1+ax\frac{B'}{B}\right)^2\right).
    \end{equation}
    Using both \eqref{eq:ident_squar} and \eqref{eq:ident_squar_diff}, we can eliminate the term $(\frac{B}{x})^{\frac{2}{a+1}}$ (which contains a non-integer power of $B$), and obtain a classical ODE on $B$. After a few simplifications, we exactly obtain \eqref{eq:diff_eq_B}. 

    We now show that any solution $B$ of \eqref{eq:diff_eq_B} with $B(0)=0$ and $B'(0)=1$ must satisfy \eqref{eq:implicit_B}, which will complete the proof. Observe that by \eqref{eq:diff_eq_B}, we necessarily have $B''(0)=0$. Moreover, the previous discussion shows that for any $\alpha>0$, the function 
\begin{equation*}
x^{-1/a}\, -\, \frac{1}{a}\int_{\bigl(x^{1/a}B(x)\bigr)^{a/(a+1)}}^{\alpha} \frac{\dd u}{u^{1+1/a}\sqrt{1+u^2}} 
\end{equation*}
is a constant $c_\alpha$. It remains to compute this constant.
We observe that
\begin{equation*}
c_{\alpha}\,=\, \frac{1}{a}\int_{\bigl(x^{1/a}B(x)\bigr)^{a/(a+1)}}^{\alpha} \PAR{1 -\frac{1}{\sqrt{1+u^2}}}\frac{\dd u}{u^{1+1/a}}\,+\, x^{-1/a}\PAR{1-\PAR{\frac{x}{B(x)}}^{\frac{1}{a+1}}}\,+\,\alpha^{-1/a}.
\end{equation*}
Since $2-\frac{1}{a}>0$, taking the limit when $x\to 0$, we deduce
\begin{equation*}
c_{\alpha}\,=\, \frac{1}{a}\int_{0}^{\alpha} \PAR{1-\frac{1}{\sqrt{1+u^2}}}\frac{\dd u}{u^{1+1/a}}\, + \, \alpha^{-1/a}\, \to\, \frac{1}{a}\int_{0}^{\infty} \PAR{1-\frac{1}{\sqrt{1+u^2}}}\frac{\dd u}{u^{1+1/a}}
\end{equation*}
as $\alpha\to \infty.$ We finally compute
\begin{eqnarray*}
\frac{1}{a}\int_{0}^{\infty} \PAR{1-\frac{1}{\sqrt{1+u^2}}}\frac{\dd u}{u^{1+1/a}} & = & \frac{1}{a}\int_0^\infty\frac{\cosh v-1}{\PAR{\sinh v}^{1+1/a}}\,\dd v\\
& = & \frac{2^{-1/a}}{a}\int_0^\infty\frac{\dd v}{\PAR{\sinh \frac{v}{2}}^{1/a-1}\PAR{\cosh \frac{v}{2}}^{1+1/a}}\\
& = & \frac{2^{-1/a}}{a}\int_0^\infty\frac{\dd x}{x^{1/2a}(1+x)^{1+1/2a}}\\
& = & \frac{2^{-1/a} \Gamma\PAR{1+ \frac{1}{a}}}{\Gamma\PAR{1+\frac{1}{2a}}}\,\Gamma\PAR{1-\frac{1}{2a}}\, =\, \frac{{\rm B} \bigl(\frac{1}{2}+\frac{1}{2a},1-\frac{1}{2a}\bigl)}{2}\, =\, \rho_a^{1/a}
\end{eqnarray*} 
where we have used the changes of variable $u=\sinh v$ and $x=\PAR{\sinh \frac{v}{2}}^2$ and the duplication formula for the Gamma function. Putting everything together shows that any solution $B$ of \eqref{eq:diff_eq_B} with $B(0)=0$ and $B'(0)=1$ satisfies \eqref{eq:implicit_B}, as required.
\endproof

We can now finish the proof of Proposition~\ref{prop:exact}. The above Lemma~\ref{prop:explicit_implicit_B} and the second equation in \eqref{eq:diff_system} show after some light algebraic computations that
\begin{equation}
\label{AB}
B(x)\,=\, x\left(\frac{G(x)}{x}\right)^{\frac{a+1}{a}}\qquad\mbox{and}\qquad A(x)\, = \,\left(\frac{G(x)}{x}\right)^{\frac{1}{a}}\sqrt{1+G(x)^2}
\end{equation}
for every $x\in (0,1/\rho_a),$ whence the required expression for the generating function $M(x).$ 
\qed

\subsection{Proof of Theorem~\ref{Mom}}
\label{subsec:proof_as}
We will use the convergent series representations
\begin{equation}
\label{x1a}
x^{-1/a} \, -\, \rho_a^{1/a} \, =\, \frac{\rho_a^{1/a} (1- x\rho_a)}{a}\, \sum_{n\geq 0} \frac{(1+1/a)_n}{(n+1)!} \, (1- x\rho_a)^n
\end{equation}
for every $x\in (0, 1/\rho_a),$ and
\begin{equation}
\label{Lagrange}
F(z) \, =\, z^{-(1+1/a)} \sum_{n\geq 0} \frac{(-1)^n (1/2)_n}{(1+a +2an) n!} \, z^{-2n}
\end{equation}
for every $z > 1.$ The first one is a simple consequence of $x^{-1/a} = \rho^{1/a} \bigl(1- (1-x\rho_a)\bigr)^{-1/a}$ and the binomial theorem, whereas the second one follows from the change of variable
\begin{equation*}
   F(z) \, =\, \frac{1}{2a}\, \int_0^{1/z^2} v^{\frac{1-a}{2a}} \frac{dv}{\sqrt{1+v^2}}\, =\, \frac{1}{2a} \sum_{n\geq 0} \frac{(-1)^n (1/2)_n}{n!}\lpa\int_0^{1/z^2} v^{\frac{1-a}{2a}+ n}\, dv\rpa,
\end{equation*}
see \eqref{eq:def_F_cv}, where in the second equality we have applied the binomial theorem to $(1+v^2)^{-1/2}$ and the involved series is easily seen to be absolutely convergent for all $z > 1$, so that the switching between the sum and the integral is justified. The Lagrange inversion formula -- see e.g.\ Appendix~E in \cite{AAR}~-- implies then after some computations the expansion
\begin{equation*}
    F^{-1}(y)\, \sim\, \left(\frac{1}{(a+1)y}\right)^{\frac{a}{a+1}} \left( 1 \,-\,\sum_{j=1}^\infty c_j\, y^{\frac{2ja}{a + 1}}\right)\qquad \mbox{as $y\to 0,$}
\end{equation*}
with 
\begin{equation*}
   c_1\, =\, \,\frac{a(a + 1)^{\frac{2a}{a + 1}}}{2(3a + 1)}\cdot
\end{equation*}
The further constants $c_p$ can be computed explicitly from \eqref{Lagrange} and, here and throughout, we have used the usual notation $\sim$ for asymptotic expansions -- see e.g.\ Appendix~C of \cite{AAR}. This entails 
\begin{equation}
\label{Complete}
F^{-1}(y)^{\frac{a+1}{a}}\, \sim\, \frac{1}{(a+1)y}\left( 1 \,-\,\sum_{j=1}^\infty {\widetilde c}_j\,  y^{\frac{2ja}{a + 1}}\right)\qquad \mbox{as $y\to 0,$}
\end{equation}
with
\begin{equation*}
   {\widetilde c}_1\,=\,\frac{(a + 1)^{\frac{3a+1}{a + 1}}}{2(3a + 1)}
\end{equation*}
and where the $\widetilde{c}_p$ can be computed recursively from the $c_p$. In order to obtain our moment estimates, we will combine Proposition \ref{prop:exact}, Equations \eqref{eq:expression_B_M_F},  \eqref{x1a} and \eqref{Complete} and the classical singularity analysis of \cite{FO}. First, observe that \eqref{Lagrange} can be expressed as a hypergeometric series, 
\begin{equation}
\label{Hypergeo}
F(z) \, =\, \frac{z^{-(1+1/a)}}{a+1}\,\pFq{2}{1}{1/2,,1/2+1/2a}{3/2+1/2a}{-1/z^2},
\end{equation}
which shows by Euler's integral formula -- see e.g.\ Theorem~2.2.1 in \cite{AAR} -- that the function $F$ is holomorphic on $\CC\setminus \bigl((-\infty,0]\cup\{\ii u, u\in[-1,1]\}\bigr)$. It will actually turn out from the alternative representation \eqref{Hypergeo2} below that this function is holomorphic on the whole cut plane $\CC\setminus (-\infty,0].$ This implies that there exists $\varepsilon > 0$ such that the function $z\mapsto F^{-1}(z)^{\frac{a+1}{a}}$ appearing in \eqref{Complete} is holomorphic on the cut open disk $\{\vert z\vert < \varepsilon\}\setminus (-\infty,0].$ Writing again 
\begin{equation*}
    z^{-1/a} - \rho_a^{1/a}\, =\, \rho_a^{1/a} \lpa\bigl(1- (1-z\rho_a)\big)^{-1/a} -1\rpa
\end{equation*}
and applying Proposition \ref{prop:exact} show that there exists $\eta > 0$ such that the function $z\mapsto M(z)$ is holomorphic on the indented open disk $\Delta_{a,\eta}\, =\,\{\vert z - 1/\rho_a\vert < \eta\}\cap \{\vert \arg(z-1/\rho_a)\vert > (1- a)\pi\}.$ Since $a > 1/2,$ we are in position to apply \cite[Thm~1]{FO} to the generating function $M$ in \eqref{eq:def_M} around the critical point $z = 1/\rho_a.$ Plugging \eqref{x1a} in \eqref{Complete} and applying Proposition \ref{prop:exact} imply
\begin{equation*}
    M(z)\, -\, \frac{2a}{(a+1)(1-z\rho_a)}\, = \, {\rm O} \bigl((1-z\rho_a)^{-\delta_a}\bigr)\qquad\mbox{as $z\to 1/\rho_a$ inside $\Delta_{a,\eta}.$}
\end{equation*}
Therefore, using $\delta_a < 1$ and applying Theorem 1 in \cite{FO}, we obtain
\begin{equation*}
   m_n\, \sim\, \frac{2a\, \rho_a^n}{a+1}
\end{equation*} 
which gives the principal estimate \eqref{Asympto1}. See Figure~\ref{fig:plot_rho} for an illustration of the above asymptotic result.

\medskip

In order to get the higher order asymptotics \eqref{Asympto}, we have to handle the odd and even moments separately. We shall begin with the odd moments and consider the renormalized sequence $\mu_n = \rho_a^{-2n-1} m_{2n+1},$ whose generating function reads
\begin{equation*}
   \sum_{n\geq 0} \mu_n u^n \, =\, \frac{B(x)}{\rho_a \,x}\, =\, \frac{1}{\rho_a}\lpa\frac{G(x)}{x}\rpa^{1+1/a}\! =\; \rho_a^{1/a} \lpa\sum_{n\geq 0} \frac{(1/2+1/2a)_n}{n!}\, (1-u)^n\rpa G(x)^{1+1/a},
\end{equation*}
where we have set $u = x^2 \rho_a^2,$ the second equality is given by \eqref{AB}, and the third equality follows from the binomial theorem applied to $x^{-1-1/a} = \rho_a^{1+1/a} \bigl((1- (1-u)\bigr)^{-(1/2+1/2a)}.$ Using \eqref{Complete} and plugging in the series representation
\begin{align}
    y\, =\,x^{-1/a} \, -\, \rho_a^{1/a} &\, =\, \frac{\rho_a^{1/a} (1- u)}{2a}\, \sum_{n\geq 0} \frac{(1+1/2a)_n}{(n+1)!} \, (1- u)^n\nonumber \\&\, =\, \frac{\rho_a^{1/a} (1-u)}{2a}\lpa 1\, +\, \frac{2a + 1}{4a}(1-u)\, +\, {\rm O} (1-u)^2\rpa
    \label{eq:expansion_y_x_1/a}
\end{align}
which is a variation on \eqref{x1a}, we deduce the expansion
\begin{eqnarray*}
\sum_{n\geq 0} \mu_n u^n & \sim & \frac{2a}{(a+1)(1-u)}\,\times\,\frac{{\displaystyle \lpa\sum_{n\geq 0} \frac{(1/2+1/2a)_n}{n!}\, (1-u)^n\rpa \lpa 1\, -\, \sum_{j=1}^\infty {\widetilde c}_j\,  y^{\frac{2ja}{a + 1}}\rpa}}{{\displaystyle \sum_{n\geq 0} \frac{(1+1/2a)_n}{(n+1)!}\, (1-u)^n}}\\
& \sim  & \frac{2a}{(a+1)(1-u)} \lpa 1 \,+\, \frac{1-u}{4a}\, +\,{\rm O}(1-u)^2 \rpa \lpa 1\, -\, \sum_{j=1}^\infty {\widetilde c}_j\,  y^{\frac{2ja}{a + 1}}\rpa\\
& \sim &  \frac{2a}{a+1}\lpa \frac{1}{1-u} \, - \, \frac{{\widehat c}_1}{(1-u)^{\delta_a}} \,+\, \frac{1}{4a}\, + \, {\widehat c}_2 (1-u)^{1-2\delta_a}\, + \, {\widehat c}_3 (1-u)^{1-\delta_a}\, +\,{\rm O}(1-u) \rpa
\end{eqnarray*}
as $u \to 1,$ with 
\begin{equation*}
   {\widehat c}_1 \, =\, \frac{{\widetilde c_1} \, \rho_a^{\frac{2}{a+1}}}{(2a)^{\frac{2a}{a+1}}}\, =\, 2^{\delta_a} \kappa_a \lpa\frac{a+1}{3a+1}\rpa
\end{equation*}
and ${\widehat c}_2, {\widehat c}_3$ are computable constants. Applying Corollary 3 in \cite{FO}, we deduce 
\begin{eqnarray*}
\mu_n & = & \frac{2a}{a+1}\lpa 1  -  \frac{{\widehat c}_1\, (\delta_a)_n}{n!}  +  \frac{{\widehat c}_2 (2\delta_a-1)_n}{n!} +  \frac{{\widehat c}_3 (\delta_a-1)_n}{n!}\, +\,{\rm O}\lpa n^{-2}\rpa\rpa  \\
& = & \frac{2a}{a+1}\lpa 1  -   \frac{{\bar c}_1 (\delta_a)_{2n+1}}{(2n+1)!} +  \frac{{\bar c}_2\, (2\delta_a)_{2n+1}}{(2n+2)!} +  \frac{{\bar c}_3\, (\delta_a)_{2n+1}}{(2n+2)!}\, +\,{\rm O} \bigl(2n+1\bigr)^{-2}\rpa 
\end{eqnarray*} 
with 
\begin{equation*}
   {\bar c}_1 \, =\, 2\kappa_a \lpa\frac{a+1}{3a+1}\rpa
\end{equation*}
and ${\bar c}_2, {\bar c}_3$ are computable constants and where in the second equality we have used \eqref{Posh}, which implies
\begin{equation*}
   \frac{(\delta_a)_n}{n!}\; =\; 2^{\frac{2a}{a+1}} \frac{(\delta_a)_{2n+1}}{(2n+1)!}\lpa 1\, +\, \frac{2a(3+a)}{(a+1)^2 n}\, +\,{\rm O}(n^{-2}) \rpa.
\end{equation*}

We now proceed to the even moments and consider the sequence $\nu_n = \rho_a^{-2n} m_{2n},$ whose generating function reads, similarly as above,
\begin{equation*}
   \sum_{n\geq 0} \nu_n u^n \, =\, A(x) \, =\,\lpa\frac{G(x)}{x}\rpa^{1/a}\sqrt{1+G(x)^2}\, =\, B(x) \lpa\sum_{n\geq 0} \frac{(-1/2)_n (-1)^n}{n!}\,G(x)^{-2n}\rpa
\end{equation*}
with the same notation $u = \rho_a^2 x^2.$ This implies the expansion
\begin{eqnarray*}
\sum_{n\geq 0} \nu_n u^n & = & \sqrt{u} \lpa\sum_{n\geq 0} \mu_n u^n\rpa \lpa\sum_{n\geq 0} \frac{(-1/2)_n (-1)^n}{n!}\,G(x)^{-2n}\rpa\\
& = & \lpa 1 \, -\, \frac{1-u}{2}\rpa\lpa \sum_{n\geq 0} \mu_n u^n\rpa   \lpa 1\, +\, \frac{1}{2G(x)^2}\, -\,\frac{1}{8G(x)^4}\rpa \, +\, {\rm O} (1-u)^2
\end{eqnarray*}
as $u \to 1.$ Using
\begin{equation*}
   \frac{1}{G(x)^2}\, =\, \PAR{(a+1) y}^{\frac{2a}{a+1}} \lpa 1 \, +\, \sum_{j\geq 1} c_j\, y^{\frac{2ja}{a+1}}\rpa^{-2}
\end{equation*}
with the expansion \eqref{eq:expansion_y_x_1/a},
and the previous estimates on the odd moments we deduce, after some simplifications,
\begin{multline*}
\sum_{n\geq 0} \nu_n u^n   \sim \\  \frac{2a}{a+1}\lpa \frac{1}{1-u} \, + \, \frac{{\widehat c}_4}{(1-u)^{\delta_a}} \,+\, \frac{1-2a}{4a} \, + \, {\widehat c}_5 (1-u)^{1-2\delta_a}\, + \, {\widehat c}_6 (1-u)^{1-\delta_a}\, +\,{\rm O}(1-u) \rpa
\end{multline*}
as $u \to 1$, with 
\begin{equation*}
   {\widehat c}_4 \, =\, \lpa\frac{a\, \rho_a^{\frac{2}{a+1}}}{3a+1}\rpa\lpa \frac{a+1}{2a}\rpa^{\frac{2a}{a+1}}\, =\, 2^{1+\delta_a} \kappa_a \lpa\frac{a}{3a+1}\rpa
\end{equation*}
and ${\widehat c}_5, {\widehat c}_6$ are computable constants. Applying again Corollary 3 in \cite{FO}, we obtain 
\begin{eqnarray*}
\nu_n  & = & \frac{2a}{a+1}\lpa 1 \, + \, \frac{{\widehat c}_4\, (\delta_a)_n}{n!} \, + \, \frac{{\widehat c}_5\, (2\delta_a-1)_n}{n!}\, + \, \frac{{\widehat c}_6\, (\delta_a-1)_n}{n!}\, +\,{\rm O}\lpa n^{-2}\rpa\rpa \\
& = & \frac{2a}{a+1}\lpa 1 \, + \,  \frac{{\bar c}_4 \, (\delta_a)_{2n}}{(2n)!}\, + \, \frac{{\bar c}_5 \,(2\delta_a)_{2n}}{(2n+1)!}\, + \, \frac{{\bar c}_6 \,(\delta_a)_{2n}}{(2n+1)!}\, +\,{\rm O}\lpa (2n)^{-2}\rpa\rpa
\end{eqnarray*} 
with 
\begin{equation*}
   {\bar c}_4 \, =\, \frac{4a\,\kappa_a}{3a+1}
\end{equation*}
and ${\bar c}_5, {\bar c}_6$ are computable constants. Putting together the expansions of odd and even moments completes the proof.
\qed

\begin{Rem} {\em Applying the transformation (2.3.12) in \cite{AAR} to the hypergeometric expression \eqref{Hypergeo} gives the series representation
\begin{equation}
\label{Hypergeo2}
F(z) \, =\, -\rho_a^{1/a}\, +\,z^{-1/a}\,\pFq{2}{1}{1/2,,-1/2a}{1-1/2a}{-z^2},
\end{equation}
which defines an analytic function on the cut disk $\{\vert z\vert \leq 1\}\setminus (-\infty,0]$. By \eqref{eq:expression_B_M_F}, this shows that the function $C(x) = B(x)^{\frac{a}{a+1}} x^{\frac{1}{a+1}}$ is the solution to the implicit equation 
\begin{equation*}
   C(x) \lpa \pFq{2}{1}{1/2,,-1/2a}{1-1/2a}{-C(x)^2}\rpa^{-a}\, =\, x,
\end{equation*}
which can be inverted by the Lagrange formula and we retrieve $B(x)\,=\, x+\frac{a+1}{2(2a-1)}x^3+\cdots.$ Observe that \eqref{Hypergeo2} also gives a direct proof of the identity
\begin{equation*}
   \rho_a^{1/a}\, =\, \frac{1}{a}\int_{0}^{\infty} \PAR{1-\frac{1}{\sqrt{1+u^2}}}\frac{\dd u}{u^{1+1/a}}\, =\, \lim_{x\downarrow 0}\lpa x^{-1/a} \, -\, F(x)\rpa,
\end{equation*}
which was used at the end of Lemma \ref{prop:explicit_implicit_B}.}
\end{Rem}

\subsection{Elementary properties of the exponential growth} We have the following property of the exponential rate of growth of the sequence $\{m_n\}$. See Figure~\ref{fig:plot_rho}.

\begin{Prop}
\label{ExpGro}
The function $a\mapsto \rho_a$ is decreasing and convex on $(1/2,1)$ from $\infty$ to 1 with
\begin{equation*}
   \rho_a\,\sim\, \frac{1}{2\sqrt{a-1/2}}\;\;\;\mbox{as $a\to 1/2$}\qquad\mbox{and}\qquad \rho_a-1\, \sim\, (1-a)\log 2\;\;\; \mbox{as $a\to 1.$}
\end{equation*}
\end{Prop}

\proof We start with the decreasing property. Setting $b = 1/2a\in (1/2,1),$ we need to show that function
\begin{equation*}
   f(b)\, =\, \lpa \frac{\Ga(1/2+b)\Ga(1-b)}{\Ga(1/2)\Ga(1)}\rpa^{\frac{1}{2b}}\, =\, \lpa \prod_{n\geq 0} \frac{(1/2+n)(1+n)}{(1/2+b+n)(1-b+n)}\rpa^{\frac{1}{2b}}
\end{equation*}
increases on $(1/2,1),$ where in the second equality we have used the classical product representation of the Gamma function -- see e.g.\ Theorem~1.1.2 in \cite{AAR}. Taking the logarithm, we are reduced to show that the function $b\mapsto b^{-1} g_n(b)$ increases on $(1/2,1)$ for each $n \geq 0,$ with 
\begin{equation*}
   g_n(b)\, =\, \log (1/2+n)\, +\, \log(1+n)\, -\, \log(1/2+b+n)\, -\, \log(1-b+n).
\end{equation*}
The function $g_n$ is strictly convex on $[1/2,1)$ with $g_n(1/2) = 0$ and
\begin{equation*}
   g_n'(b)\, =\, \frac{2b-1/2}{(1/2+b+n)(1-b+n)}\, >\, 0\qquad \mbox{for all $b\in[1/2,1).$}
\end{equation*}
Therefore, the function
\begin{equation*}
   \frac{g_n(b)}{b}\, =\, \lpa\frac{b-1/2}{b}\rpa\times\lpa\frac{g_n(b) -g_n(1/2)}{b-1/2}\rpa
\end{equation*}
increases on $(1/2,1)$ as the product of two positive increasing functions. Moreover, the function $a\mapsto 2a g_n(1/2a)$ has second derivative $(1/2a^3) g_n''(1/a) > 0$ on $(1/2,1]$ and is hence convex. This implies that the mapping $a\mapsto \rho_a$ is log-convex and hence convex on $(1/2,1).$ Last, we compute
\begin{equation*}
   (a-1/2)^a \rho_a\, =\, \lpa \frac{a}{\sqrt{\pi}} \, \Gamma\lpa \frac{1}{2} + \frac{1}{2a}\rpa \Gamma \lpa 2 - \frac{1}{2a}\rpa \rpa^a\, \to\, \frac{1}{2}\quad\mbox{as $a\to 1/2,$}
\end{equation*}
whence the first asymptotics. Setting $\psi(z)=\frac{\Gamma'(z)}{\Gamma(z)}$ for the usual digamma function, we have 
\begin{equation*}
   \lpa \log \rho_a\rpa'\, = \, \frac{1}{2a}\bigl( \log\rho_a \, +\, \psi(1-1/2a)\, -\, \psi(1/2+1/2a)\bigr)\, \to \, \frac{\psi(1/2)-\psi(1)}{2}\, =\, -\log 2\quad \mbox{as $a\to 1$}
\end{equation*}
whence the second asymptotics.

\qed

\begin{Rem} {\em The decreasing character of $a \mapsto \rho_a$ on $(1/2,\infty)$ also follows from Hölder's inequality and the following representation:
\begin{equation*}
   \rho_a \, =\, \esp \lcr\lpa \frac{\GG_{1/2}}{\GG_1}\rpa^{\!\frac{1}{2a}}\rcr^a
\end{equation*}
where $\GG_t$ stands for the standard Gamma random variable with parameter $t$ and the quotient inside the expectation is independent.}

\medskip

{\em As a second note, applying Theorem 1.6.2 (ii) in \cite{AAR} it is possible to show that 
\begin{equation*}
   \log \rho_a\, = \, -\log 2 \, +\, \int_0^\infty e^{-ax} \,f(x)\, dx
\end{equation*}
for all $a\in (1/2,\infty),$ where  
\begin{eqnarray*}
f(x) & = & \frac{1}{x^2}\sum_{n\geq 0}\lpa 2 - e^{-\frac{x}{2n+1}} \lpa 1 + \frac{x}{2n+1}\rpa \, +\, e^{\frac{x}{2n+2}}\lpa \frac{x}{2n+2}-1\rpa\rpa
\end{eqnarray*}
is easily shown to be positive on $(0,\infty).$ This implies that the mapping $a\mapsto \rho_a$ is logarithmically completely monotone and hence completely monotone on $(1/2,\infty),$ with limit $1/e^2$ at infinity. 
}
\end{Rem}

\subsection{Density asymptotics at the logarithmic scale} We end this section with a logarithmic estimate of $\varphi(x)$ as $x\to\pm\infty.$ This less precise version of Theorems \ref{Positano} and \ref{Negus} can be quickly derived from Theorem \ref{Mom} and Kasahara's Tauberian theorem for densities. 

\begin{Prop}
\label{Log}
As $x\to\infty$, one has
\begin{equation*}
   \log \varphi(x)\, \sim\, \log\varphi(-x)\, \sim\, -(1-a)  \lpa \frac{a^a x}{\rho_a}\rpa^{\frac{1}{1-a}}.
\end{equation*}
\end{Prop}

\proof
Along the proof we will use the classical Mittag-Leffler function $E_a$ defined by
\begin{equation*}
    E_a(z)\, =\, \sum_{n\geq 0} \frac{z^n}{\Ga(1+an)},
\end{equation*}
and Mittag-Leffler random variables $M_a$ with moment generating function
$\esp [e^{rM_a}] = E_a(r)$.

We begin the proof of Proposition~\ref{Log} with the estimate on the positive axis. Introduce the moment generating function $\Psi(r)= \esp[e^{rL_1}]=\sum_{n\geq 0} \bigl(\frac{\esp[L_1^n]}{n!}\bigr) r^n$ as in \eqref{eq:def_MGF}, 
which is analytic on $\CC$ by Stirling's formula combined with the equivalence
\begin{equation*}
   \frac{\esp[L_1^n]}{n!}\,\sim\, \frac{2a\, \rho_a^n}{(a+1)\,\Ga(1+an)}
\end{equation*}
given by Theorem \ref{Mom}. The latter estimate also implies
\begin{equation*}
   \Psi (r) \,\sim\, \frac{2a}{a+1}\, \sum_{n\geq 0} \frac{(\rho_a r)^n}{\Ga(1+an)}\, =\, \frac{2a\, E_a(\rho_a r)}{a+1}\,\sim\, \frac{2\, e^{(\rho_a r)^{1/a}}}{a+1}\qquad\mbox{as $r\to\infty,$}
\end{equation*}
where in the second equivalence we have used the well-known estimate for the Mittag-Leffler function $E_a$ on the positive half-line, to be found e.g.\ in Formula~(3.4.14) of \cite{GKMR}. Combining this estimate on $\Psi(r)$ with Kasahara's theorem for densities -- see Theorem 4.12.11 in \cite{BGT} -- yields after some easy algebra the required estimate  
\begin{equation*}
   \log \varphi(x)\,\sim\, -(1-a)  \lpa \frac{a^a x}{\rho_a}\rpa^{\frac{1}{1-a}}\qquad \mbox{as $x\to\infty.$}
\end{equation*}

We now proceed to the estimate on the negative axis. Let $M_a$ be a Mittag-Leffler random variable and let $B_a$ be an independent Bernoulli random variable with parameter $\frac{2a}{a+1}\cdot$ The random variable $Y_a = (\rho_a B_a)\times M_a$ has positive integer moments 
\begin{equation*}
   \esp [Y^n_a] \, =\, \frac{2a \,\rho_a^n\,n!}{(1+a)\,\Gamma(1+an)}\cdot
\end{equation*}
Consider the function
\begin{equation*}
   {\widetilde \Psi}(r)\, =\, \Psi(-r) \, -\, \esp[e^{-r Y_a}] \, =\, \sum_{n\geq 0} \frac{\mu_n}{n!}\, r^n
\end{equation*}
with
\begin{eqnarray}
\mu_n & = & (-1)^n \lpa \esp[L_1^n] - \frac{2a \,\rho_a^n\,n!}{(1+a)\,\Gamma(1+an)}\rpa\nonumber \\ 
& = & \frac{2a \kappa_a\, \rho_a^n\,(\delta_a)_n}{(a+1)\,\Gamma(1+an)}\lpa 1\, +\, (-1)^n\frac{a-1}{3a+1}\, +\, {\rm O}\lpa n^{-\delta_a}\rpa\rpa\, \asymp\, \frac{\rho_a^n\,(\delta_a)_n}{\Gamma(1+an)}\label{eq:bounds_Prab}
\end{eqnarray}
by Theorem \ref{Mom} and \eqref{Posh}. Here and throughout, for two real sequences $\{u_n\}$ and $\{v_n\},$ the notation $u_n \asymp v_n$ means that there exist an index $n_0\geq 1$ and two constants $c, C > 0$ such that $cv_n \leq u_n \leq C u_n$ for all $n\geq n_0$. The upper and lower bounds in \eqref{eq:bounds_Prab} follow from $0 < \frac{1-a}{3a+1} < 1.$ 

Introduce the following two-parameter generalization of the Mittag-Leffler function, called the Prabhakar function
\begin{equation}
\label{Prabhakar}    
E^\gamma_{\a,\b} (z)\, =\, \sum_{n\geq 0} \frac{(\gamma)_n\, z^n}{n!\, \Gamma(\b +\a n)},
\end{equation}
which is analytic on $\CC$ for all $\a,\b,\gamma > 0$. Using \eqref{eq:bounds_Prab}, we deduce that
\begin{equation*}
   \log \esp \lcr e^{r\vert L_1\vert} \Un_{\{L_1 \leq 0\}}\rcr \, \sim\,\log {\widetilde \Psi} (r)\, \sim\, \log E^{\delta_a}_{a,1} (\rho_a r)\, \sim\, (\rho_a r)^{1/a}\, \to\, \infty\quad \mbox{as $r \to \infty,$}
\end{equation*}
where the first equivalence is clear by the positivity of $Y_a$ and for the third equivalence, we have used the global estimate presented in Theorem~3 of \cite{GG} -- see also \cite{W40}. As in the case of the positive axis, the conclusion follows from Kasahara's Tauberian theorem for densities.
\endproof

\begin{Rem}{\em If the sequence $\{m_n\}$ were a Hamburger moment sequence, then we would have an independent factorization 
$$L_1\, \elaw \, X\times M_a$$
where $X$ is the real random variable having integer moments $m_n,$ whose asymptotic behaviour would show that the right end of the support of $X$ equals $\rho_a$ and the left end equals $-\rho_a.$ The statement of Proposition \ref{Log} would then follow in a more elementary way by convolution from the known asymptotics
$$f_{M_a} (x)\, \sim\, \kappa_a x^{\frac{2a-1}{2(1-a)}} e^{-(1-a) (a^a x)^{\frac{1}{1-a}}}\qquad \mbox{as $x\to\infty,$}$$
where $f_{M_a}$ is the density function of $M_a$ -- combine e.g.\ Formula (14.35) and Exercise 29.18 in \cite{Sato} for the latter estimate -- and $\kappa_a > 0$ some constant. Unfortunately, simulations show that $\{m_n\}$ is not a Hamburger moment sequence.}
\end{Rem} 

\section{Proof of Theorem \ref{Positano}}

\subsection*{Strategy of the proof}
In this section we prove Theorem \ref{Positano}, with the following explicit value of the prefactor $c_a$:
\begin{equation}
\label{eq:value_c_a}
   c_a \, =\, \sqrt{\frac{2}{\pi(1-a^2)(1+a)}} \lpa \frac{a}{\rho_a}\rpa^{\frac{1}{2(1-a)}}.
\end{equation}
Introducing the random variable $X = -L_1/\rho_a$ with density function ${\widetilde \varphi} (x) = \rho_a \, \varphi (-\rho_a x)$ on $\rl,$ the required estimate amounts to
\begin{equation}
\label{Posit}
{\widetilde \varphi}(-x)\, \sim\, \sqrt{\frac{2}{\pi(1-a^2)(1+a)}}\; a^{\frac{1}{2(1-a)}} \, x^{\frac{2a-1}{2(1-a)}}\, e^{-(1-a)  \lpa a^a x\rpa^{\frac{1}{1-a}}}, \qquad x\to\infty.
\end{equation}
To prove \eqref{Posit}, we will use the strong Tauberian theorem given in \cite[Thm~3]{FY}. This result allows us to derive the asymptotics of the density of a random variable $X$, after checking three technical conditions involving the following functionals of the distribution of $X$:
\begin{equation}
\label{eq:def_omega_eta}
\omega(r)\, =\, \esp[e^{-r X}], \quad \xi(r)\, =\, -\frac{\omega'(r)}{\omega(r)}\quad\mbox{and}\quad \eta^2(r)\, =\, -\xi'(r)\, =\, \frac{\omega''(r)}{\omega(r)}\, -\, \xi^2(r).
\end{equation}
We divide the proof into two steps: first, we will derive the asymptotic behaviour of the functions introduced in \eqref{eq:def_omega_eta}. In the second step, we will use these estimates to verify the three conditions of \cite[Thm~3]{FY}.

\subsection*{Detailed proof of Theorem \ref{Positano}}
\begin{proof}[First step: asymptotics of the functions \eqref{eq:def_omega_eta}]
Let us first note that the functions defined in \eqref{eq:def_omega_eta}
 are all smooth real functions, since $\omega(r)$ is smooth and never vanishes on $\rl$ as a moment generating function with infinite radius of convergence. As a consequence of \eqref{Asympto} in Theorem~\ref{Mom} and the estimate $an\Ga (1+ an) \sim\Ga(2+an),$ we have the expansion
\begin{equation*}
   \frac{\esp[X^n]}{n!}\, = \,\frac{2a\, (-1)^n}{(a+1)}\!\lpa \frac{1}{\Ga(1+an)} + \kappa_a \frac{(\delta_a)_n ( (-1)^n  - \frac{1-a}{3a+1})}{n!\,\Ga(1+an)}\, +\, {\rm O}\lpa \frac{(2\delta_a)_n}{n!\,\Ga(2+an)}\!\rpa\!\rpa
\end{equation*}
as $n\to\infty.$ Recalling the above notation \eqref{Prabhakar} of the Prabhakar function $E^\gamma_{\a,\b} (z),$ we deduce
\begin{eqnarray}
\label{Base}
\omega (r) & = & \sum_{n\geq 0}\lpa\frac{(-1)^{n}\, \esp[X^n]}{n!}\rpa r^n\nonumber\\
& = & \frac{2a}{a+1}\lpa E_a(r)  + \kappa_a \lpa E_{a,1}^{\delta_a} (-r) + \frac{a-1}{3a+1}\, E_{a,1}^{\delta_a} (r) \rpa \rpa +  {\rm O}\lpa E_{a,2}^{2\delta_a} (r) \rpa\nonumber\\
& = & \frac{2\, e^{r^{1/a}}}{a+1}\lpa 1 \, -\,{\widetilde \kappa}_a r^{-2/(a+1)}\, +\, {\rm O}\lpa r^{-4/(a+1)}\rpa \rpa,
\end{eqnarray}
with the notation 
\begin{equation*}
   {\widetilde \kappa}_a \, =\, \frac{\kappa_a a^{1-\delta_a}(1-a)}{3a+1}
\end{equation*}
and where for the equality \eqref{Base} we have used the estimates
\begin{equation*}
    E_{\a,\beta}^{\gamma} (r)\, =\, \frac{r^{\frac{\gamma -\beta}{\a}}}{\a^\gamma} \, e^{r^{1/\a}}\lpa 1\, +\, {\rm O} (r^{-1/\a})\rpa\qquad\mbox{and}\qquad E_{\a,\b}^{\gamma} (-r) \, = \, {\rm O} (r^{-\gamma})\qquad\mbox{as $r\to\infty,$}
\end{equation*}  
which are valid for every $\a,\b,\gamma > 0$ -- see Theorem~3 in \cite{GG}. A similar argument yields 
\begin{equation*}
    \frac{\esp[X^{n+1}]}{n!}\, = \, \frac{2\, (-1)^{n+1}}{(a+1)}\lpa  \frac{1}{\Ga(a+an)}\, +\, \kappa_a \frac{(\delta_a)_n \bigl( (-1)^{n+1}  - \frac{1-a}{3a+1}\bigr)}{n!\,\Ga(a+an)} \, + \,{\rm O}\lpa \frac{(2\delta_a)_n}{n!\,\Ga(1+a+an)}\rpa\rpa
\end{equation*}
which implies, by the same estimates on the Prabhakar function, 
\begin{eqnarray}
\label{Base1}
\omega' (r) & = & \sum_{n\geq 0}\lpa\frac{(-1)^{n+1}\, \esp[X^{n+1}]}{n!}\rpa r^n\nonumber\\
& = & \frac{2a}{a+1}\lpa E_{a,a}^1(r)  + \kappa_a \lpa E_{a,a}^{\delta_a} (-r) + \frac{a-1}{3a+1}\, E_{a,a}^{\delta_a} (r) \rpa\rpa  +  {\rm O}\lpa E_{a,1+a}^{2\delta_a} (r) \rpa\nonumber\\
& = & \frac{2\, r^{1/a-1}\, e^{r^{1/a}}}{a(a+1)}\lpa 1 \, -\,{\widetilde \kappa}_a r^{-2/(a+1)}\, +\, {\rm O}\lpa r^{-4/(a+1)}\rpa \rpa.
\end{eqnarray}
Combining \eqref{Base} and \eqref{Base1} gives the first estimate
\begin{equation}
\label{XiEst}
\xi(r)\, =\, -\frac{r^{1/a-1}}{a}\lpa 1 \, + \, {\rm O}\lpa r^{-4/(a+1)}\rpa\rpa.
\end{equation}
Then, after some simplifications, the details of which are left to the reader, we obtain the estimate 
\begin{equation*}
\frac{\esp[X^{n+2}]}{n!}\, = \, \frac{2\, (-1)^n}{a(a+1)}\,\lpa \frac{1}{\Ga(2a-1+an)}\, +\, \kappa_a\, \frac{(\delta_a)_n \bigl((-1)^n\, -\,\frac{1-a}{3a +1}\bigr)}{n!\,\Ga(2a-1+an)}
\, +\, {\rm O}\lpa\frac{(2\delta_a)_n}{n!\,\Ga(2a+an)}\rpa\rpa,
\end{equation*}
which yields similarly 
\begin{equation*}
\omega''(r)\, =\, \sum_{n\geq 0}\lpa\frac{(-1)^n \,\esp[X^{n+2}]}{n!}\rpa r^n \, =\, \frac{2\, r^{2/a-2}\, e^{r^{1/a}}}{a^2(a+1)}\lpa 1 \, -\,{\widetilde \kappa}_a r^{-2/(a+1)}\, + \, {\rm O}\lpa r^{-4/(a+1)}\rpa\rpa.
\end{equation*}
Putting this estimate together with \eqref{Base} and \eqref{XiEst} implies
\begin{equation*}
   \eta^2(r)\, =\, \frac{r^{2/a-2}}{a^2}\lpa \frac{ 3a\, r^{-1/a}}{2}\, + \, {\rm O}\lpa r^{-4/(a+1)}\rpa\rpa,
\end{equation*}
whence our final estimate
\begin{equation*}
    \eta(r)\, =\, \frac{\sqrt{1-a}\, r^{1/2a-1}}{a}\lpa 1 \, + \, {\rm O}\lpa r^{-4/(a+1)}\rpa\rpa,
\end{equation*}
concluding the first part of the proof.
\end{proof}

\begin{proof}[Second step: checking the assumptions of \cite{FY}]
We are now in position to apply the strong Tauberian theorem of \cite{FY}. Specifically, with the notation of \cite{FY}, if $F$ stands for the distribution function of the above random variable $X = -L_1/\rho_a$,  one has $R=-\infty$ and $S = \infty$ so that $\omega (r)$ satisfies (2) therein. The above asymptotics of $\eta(r)$ clearly implies the slow variation condition B in \cite{FY}, and some algebra combined with the two above equivalents for $\omega(r)$ and $\eta(r)$ shows that our required \eqref{Posit} amounts to (14) therein. Hence, by Theorem 3 in \cite{FY}, our proof will be complete if the domination condition (11) therein is fulfilled, and we need to establish the existence of some $g\in\cL_1(\rl)$ such that 
\begin{equation*}
   \lpa \frac{r\eta(r)}{\sqrt{r^2\eta^2(r) + s^2}}\rpa \lva \frac{\omega(r + \ii s/\eta(r))}{\omega (r)}\rva \, \leq\, g(s)
\end{equation*}
for all $s\in\rl$ and $r$ large enough.

The estimate \eqref{Base} remains valid in the half-plane $\{\Re(z) \geq 0\},$ and applying Theorem 3 in \cite{GG} to the functions $E_a(r + \ii s/\eta(r))$, $E_{a,1}^{\delta_a} (r + \ii s/\eta(r))$, $E_{a,1}^{\delta_a} (-r - \ii s/\eta(r))$ 
and $E_{a,1}^{2\delta_a}(r + \ii s/\eta(r)),$ we obtain
\begin{equation*}
   \lva\omega (r+ \ii s/\eta(r))\rva \, \leq\, \lpa \frac{2\, \lva e^{(r+ \ii s/\eta(r))^{1/a}}\rva}{a+1}\, +\, \frac{1}{\Ga(1- a\delta_a)}\lpa \frac{\eta(r)}{\sqrt{r^2\eta^2(r) + s^2}}\rpa^{\delta_a}\rpa \lpa 1 \, + \, {\rm O}(r^{-2a/(a+1)})\rpa 
\end{equation*}
where the Landau symbol ${\rm O}$ does not depend on $s$ since $\vert r+ \ii s/\eta(r)\vert \geq r$ for all $r\geq 0$ and $s\in\rl.$ We hence need to show that there exists $g\in\cL_1(\rl)$
such that 
\begin{equation*}
   \lpa \frac{r\eta(r)}{\sqrt{r^2\eta^2(r) + s^2}}\rpa\lpa \lva e^{(r+ \ii s/\eta(r))^{1/a} -\, r^{1/a}}\rva \, + \, e^{-r^{1/a}}\lpa \frac{\eta(r)}{\sqrt{r^2\eta^2(r) + s^2}}\rpa^{\delta_a} \rpa\,\leq\, g(s)
\end{equation*}
uniformly in $s\in\rl$ for all $r$ large enough. Since $r\eta(r)\to \infty$ and $r e^{-r^{1/a}} \eta(r)^{1+\delta_a}\to 0$ as $r\to\infty,$ one has 
\begin{equation*}
   r e^{-r^{1/a}}\lpa \frac{\eta(r)}{\sqrt{r^2\eta^2(r) + s^2}}\rpa^{1+\delta_a}\,\leq\, \lpa 1+s^2\rpa^{-(1+\delta_a)/2}\, \in\,\cL_1(\rl)
\end{equation*}
uniformly in $s\in\rl$ for all $r$ large enough, and we are reduced to show that there exists $r_0 > 0$ and $h\in\cL_1(\rl)$ such that 
\begin{equation}
\label{Hyper}
\lva e^{(r+ \ii s/\eta(r))^{1/a} -\, r^{1/a}}\rva\,\leq\, h(s)
\end{equation}
for all $s\in\rl$ and $r\geq r_0.$ Setting $z = z(s,r) = s^2/r^2\eta^2(r)\geq 0$ for concision, the binomial theorem implies
\begin{eqnarray*}
\log \lva e^{(r+ \ii s/\eta(r))^{1/a} -\, r^{1/a}}\rva & = & r^{1/a} \sum_{n\geq 1} \frac{(-1/a)_{2n}}{(2n)!}\, (-z)^n\\
& = & - \lpa\frac{(1-a)z\, r^{1/a}}{a^2}\rpa \sum_{n\geq 0} \frac{(2-1/a)_{2n}}{(2n+2)!}\, (-z)^n\\
& = & - \lpa\frac{(1-a)z\, r^{1/a}}{2 a^2}\rpa \sum_{n\geq 0} \frac{(1-1/2a)_n (3/2 -1/2a)_n}{(3/2)_n (n+1)!}\, (-z)^n 
\end{eqnarray*}
with the usual notation for the Pochhammer symbol $(x)_n$ and where in the last equality we have used the Gauss multiplication formula given e.g.\ at the beginning of Chapter 1.5 in \cite{AAR}. Introducing the hypergeometric function 
\begin{equation*}
   F_a(u)\, =\, \pFq{2}{1}{1-1/2a,,3/2-1/2a}{3/2}{u}\, =\, \sum_{n\geq 0} \frac{(1-1/2a)_n (3/2 -1/2a)_n}{(3/2)_n n!}\, u^n,
\end{equation*}
this function is analytic in $\CC$  cut along the real axis from $1$ to $\infty$ since $3/2 > 3/2 - 1/2a > 0$, see e.g.\ Theorem~2.2.1 in \cite{AAR}, and we have
\begin{equation*}
   \sum_{n\geq 0} \frac{(1-1/2a)_n (3/2 -1/2a)_n}{(3/2)_n (n+1)!}\, (-z)^n\, =\, \frac{1}{z} \int_0^z F_a(-u)\, du.
\end{equation*}
Moreover, Pfaff's formula -- see e.g.\ Theorem~2.2.5 in \cite{AAR} -- and the fact that $1-1/2a > 0$ imply
\begin{equation*}
   F_a(-u) \, =\, (1+u)^{1/2a -1} \pFq{2}{1}{1-1/2a,1/2a}{3/2}{\frac{u}{u+1}}\, > \, 0
\end{equation*} 
and
\begin{eqnarray*}
-F_a'(-u) & = & - \lpa \frac{(2a-1)(3a -1)}{3a^2}\rpa \pFq{2}{1}{2-1/2a,,5/2-1/2a}{5/2}{-u} \\
& = & - \lpa \frac{(2a-1)(3a -1)(1+u)^{1/2a -2}}{3a^2}\rpa  \pFq{2}{1}{2-1/2a,1/2a}{5/2}{\frac{u}{u+1}}\, < \, 0
\end{eqnarray*}
for all $u\in\rl.$  Since $1/r^2\eta^2 (r) = {\rm O}(r^{-1/a}),$ there exists $r_1 > 0$ such that $z = z(s,r) < s^2$ for all $s\in\rl$ and $r\geq r_1,$ and by concavity we then obtain
\begin{equation*}
   \frac{1}{z} \int_0^z F_a(-u)\, du\, \geq\, \frac{1}{s^2} \int_0^{s^2} F_a(-u)\, du.
\end{equation*}
Putting everything together, we see that there exists $r_0 \geq r_1$ such that 
\begin{equation*}
   \log \lva e^{(r+ \ii s/\eta(r))^{1/a} -\, r^{1/a}}\rva \,\leq\, - \lpa\frac{(1-a)\, r^{1/a-2}}{2 a^2\,\eta^2(r)}\rpa \int_0^{s^2} F_a(-u)\, du\, \leq\, -\frac{1}{4}\, \int_0^{s^2} F_a(-u)\, du
\end{equation*}
for every $s\in\rl$ and $r\geq r_0,$ where the second inequality follows from the above equivalent for $\eta(r).$ The aforementioned Pfaff formula combined with the Gauss summation formula given e.g.\ in Theorem~2.2.2 of \cite{AAR} show that 
\begin{equation*}
   F_a(-u)\,\sim\, \lpa \frac{\pi}{2 \Ga (1+ 1/2a) \Ga(3/2- 1/2a)} \rpa u^{1/2a -1}
\end{equation*}
as $u\to\infty,$ and integrating this equivalent finally shows that there exists $c > 0$ such that
\begin{equation*}
   \log \lva e^{(r+ \ii s/\eta(r))^{1/a} -\, r^{1/a}}\rva \,\leq\, - c \vert s\vert^{1/a}
\end{equation*}
for all $s\in\rl$ and $r\geq r_0,$ which gives \eqref{Hyper} and completes the proof.
\end{proof}

\begin{Rem}{\em In the case $a\in (1/2,a_0],$ a shorter proof of Theorem \ref{Positano} can be provided using Corollary \ref{LC2} and Theorem 2 in \cite{FY}.}
\end{Rem}

\section{Proof of Theorem \ref{Negus}} 
In this section we prove Theorem \ref{Negus}, with the following explicit value of the prefactor $\widehat c_a$:
\begin{equation}
\label{eq:value_c_a_hat}
   {\widehat c}_a \, =\, \frac{1}{\sqrt{2\pi(1-a)}} \lpa \frac{\rho_a}{a}\rpa^{\frac{3a-1}{2(1-a^2)}}\frac{\rho_a^{\frac{2}{a+1}}}{2\,(a+1)^{\frac{1-a}{1+a}}\, \Gamma\lpa\frac{1-a}{1+a}\rpa}\cdot
\end{equation}
The proof is analogous to Theorem \ref{Positano} except that here we need an estimate of $\xi(-r)$ and $\eta(-r)$ as $r\to\infty,$ with the previous notation \eqref{eq:def_omega_eta}. Using \eqref{Asympto} up to the fifth order and reasoning as in the above proof, we first obtain
\begin{eqnarray*}
\omega (-r) & = & \frac{2a}{a+1}\lpa \kappa_a\, E_{a,1}^{\delta_a} (r)
\, +\, a\kappa_1^-  E_{a,2}^{2\delta_a} (r)\, + \, a\kappa_2^- E_{a,2}^{\delta_a}(r) \,  + \, {\rm O}\lpa E_{a,3}^{1} (r) \rpa\rpa \\
& = & \frac{2\, a^{\frac{2a}{a+1}} \kappa_a \,r^{-2/(a+1)}\, e^{r^{1/a}}}{a+1}\lpa 1 
\, +\, {\bar \kappa_1} r^{-2/(a+1)}\, + \, {\bar \kappa_2} r^{-1/a} \, +\, {\rm O} \lpa r^{-2/a(a+1)}\rpa \rpa
\end{eqnarray*}
with
\begin{equation*}
   {\bar \kappa_1} \, =\, \frac{a^{\frac{2a}{a+1}}\kappa_1^-}{\kappa_a}\qquad\mbox{and}\qquad {\bar \kappa_2} \, =\, \frac{a \kappa_2^-}{\kappa_a}\cdot
\end{equation*}
Similarly, using the expansion
\begin{multline*}
\frac{\esp[X^{n+1}]}{n!}\, = \, \frac{2\, (-1)^{n+1}}{(a+1)}\lpa  \frac{1}{\Ga(a+an)}\, +\, \kappa_a \frac{(\delta_a)_n \bigl( (-1)^{n+1}  - \frac{1-a}{3a+1}\bigr)}{n!\,\Ga(a+an)}
\right.\\ \left.\, +\, \frac{a\,(2\delta_a)_n\,\bigl((-1)^{n+1} \kappa_1^- + \kappa_1^+\bigr)}{n!\,\Ga(1+a+an)}\, + \, \frac{a\,(\delta_a)_n \bigl((-1)^{n+1} \kappa_2^- + \kappa_2^+-\bigr)}{n!\,\Ga(1+a+an)}\right. \\
\left.\, - \,\frac{2a^2\,\kappa_a\, (\delta_a)_n \bigl( (-1)^{n+1}  - \frac{1-a}{3a+1}\bigr)}{(a+1)\,n!\,\Ga(1+a+an) }\, + \,{\rm O}\lpa\frac{1}{\Ga(2+a+an)} \rpa\rpa,
\end{multline*}
we obtain
\begin{eqnarray*}
\omega' (-r) & = & \frac{2}{a+1}\lpa \kappa_a\, E_{a,a}^{\delta_a} (r)
\, +\, a\kappa_1^-  E_{a,1+a}^{2\delta_a} (r)\, + \, a\lpa \kappa_2^- - \frac{2a\,\kappa_a}{a+1}\rpa E_{a,1+a}^{\delta_a}(r) \,  + \, {\rm O}\lpa E_{a,2+a}^{1} (r) \rpa\rpa \\
& = & \frac{2\, \kappa_a\, r^{1/a-1-2/(a+1)}\, e^{r^{1/a}}}{a^{\delta_a}(a+1)}\lpa 1 
\, +\, {\bar \kappa_1} r^{-2/(a+1)}\, + \, {\bar \kappa_3} r^{-1/a} \, +\, {\rm O} \lpa r^{-2/a(a+1)}\rpa \rpa
\end{eqnarray*}
with
\begin{equation*}
   {\bar \kappa_3} \, =\, {\bar \kappa_2} \, -\, \frac{2a^2}{a+1}\cdot
\end{equation*}
This leads to
\begin{equation}
\label{XiEstM}
\xi(-r)\, =\, -\frac{r^{1/a-1}}{a}\lpa 1 \, -\, \frac{2a^2 \, r^{-1/a}}{a+1}\, + \, {\rm O}\lpa r^{-2/a(a+1)}\rpa\rpa.
\end{equation}
We next compute
\begin{multline*}
\frac{\esp[X^{n+2}]}{n!}\, = \, \, = \, \frac{2\, (-1)^n}{a(a+1)}\,\lpa \frac{1}{\Ga(2a-1+an)}\, +\, \kappa_a\, \frac{(\delta_a)_n \bigl((-1)^n\, -\,\frac{1-a}{3a +1}\bigr)}{n!\,\Ga(2a-1+an)}
\right.\\ \left.\, +\, \frac{a\,(2\delta_a)_n\,\bigl((-1)^{n} \kappa_1^- + \kappa_1^+\bigr)}{n!\,\Ga(2a+an)}\, + \, \frac{a\,(\delta_a)_n \bigl((-1)^{n} \kappa_2^- + \kappa_2^+-\bigr)}{n!\,\Ga(2a+an)}\right. \\
\left.\, + \,\frac{(1-5a^2)\,\kappa_a\, (\delta_a)_n ( \bigl(-1)^{n}  - \frac{1-a}{3a+1}\bigr)}{(a+1)\,n!\,\Ga(2a+an) }\, + \,{\rm O}\lpa\frac{1}{\Ga(1+2a+an)} \rpa\rpa,
\end{multline*}
which implies
\begin{eqnarray*}
\omega'' (-r) & = & \frac{2}{a(a+1)}\lpa \kappa_a\, E_{a,2a-1}^{\delta_a} (r)
\, +\, a\kappa_1^-  E_{a,2a}^{2\delta_a} (r)\right.\\
& & \qquad\qquad\qquad\qquad \left.\, + \, \lpa a\kappa_2^- + \frac{(1-5a^2)\,\kappa_a}{a+1}\rpa E_{a,2a}^{\delta_a}(r) \,  + \, {\rm O}\lpa E_{a,1+2a}^{1} (r) \rpa\rpa \\
& = & \frac{2\, \kappa_a\, r^{2/a-2-2/(a+1)}\, e^{r^{1/a}}}{a^{1+\delta_a}(a+1)}\lpa 1 
\, +\, {\bar \kappa_1} r^{-2/(a+1)}\, + \, {\bar \kappa_4} r^{-1/a} \, +\, {\rm O} \lpa r^{-2/a(a+1)}\rpa \rpa
\end{eqnarray*}
with
\begin{equation*}
   {\bar \kappa_4} \, =\, {\bar \kappa_2} \, +\, \frac{1-5a^2}{a+1}\cdot
\end{equation*} 
This leads to
\begin{equation*}
   \frac{\omega''(-r)}{\omega(-r)}\, =\, \frac{r^{2/a-2}}{a^2}\lpa 1 \, +\, \frac{(1-5a^2) \, r^{-1/a}}{a+1}\, + \, {\rm O}\lpa r^{-2/a(a+1)}\rpa\rpa
\end{equation*}
and finally, combining this estimate with \eqref{XiEstM}, to
\begin{equation}
\label{Nagast}
\eta(-r)\, =\, \sqrt{\frac{\omega''(-r)}{\omega(-r)}\, -\, \xi^2(-r)}\, =\, \frac{\sqrt{1-a}\, r^{1/2a-1}}{a}\lpa 1 \, + \, {\rm O}\lpa r^{-(1-a)/a(a+1)}\rpa\rpa.
\end{equation}
Similarly as on the positive axis, we can then apply Theorem 3 in \cite{FY} and obtain, with the previous notation for the function ${\widetilde \varphi},$ the estimate 
\begin{equation*}
   {\widetilde \varphi}(x)\, \sim\, \lpa \frac{a^{\frac{1-3a}{2(1-a^2)}}\rho_a^{\frac{2}{a+1}}}{\sqrt{8\pi(1-a)} (a+1)^{\delta_a} \Gamma\lpa \frac{1-a}{1+a}\rpa}\rpa x^{\frac{2a^2-3a-1}{2(1-a^2)}}\, e^{-(1-a)  \lpa a^a x\rpa^{\frac{1}{1-a}}}, \qquad x\to\infty,
\end{equation*}
which completes the proof.
\qed
  
\section*{Acknowledgement}
The research of HG was funded by Natural Sciences and Engineering Research Council of Canada, discovery grant RGPIN-2020-07239. KR is partially supported by the RAWABRANCH project ANR-23-CE40-0008, funded by the French National Research Agency, and by Centre Henri Lebesgue, programme ANR-11-LABX-0020-01. The authors thank Tony Guttmann for his inspiring numerical computations on the model. The authors would like to thank the Isaac Newton Institute for Mathematical Sciences, Cambridge, for support and hospitality during the programme ``Stochastic systems for anomalous diffusion''. This work was supported by EPSRC grant no EP/R014604/1.

\bibliographystyle{siam}
\bibliography{biblio}

\end{document}